\documentclass[12pt]{article}
\usepackage{latexsym}
\usepackage{amsfonts}
\usepackage{amsmath}
\usepackage{amssymb}
\usepackage{graphicx}
\usepackage[all]{xy}
\newtheorem{thm}{Theorem}

\newtheorem{lemma}{Lemma}

\newtheorem{prop}{Proposition}

\newtheorem{remark}{Remark}

\newenvironment{proof}{\medskip \noindent
{\bf Proof.}}{\hfill \rule{.5em}{1em}
\\}

\def\bea{\begin{eqnarray*}}
\def\eea{\end{eqnarray*}}
\def\be{\begin{equation}}
\def\ee{\end{equation}}

\def\k{K\"{a}hler  }

\begin{document}

\title{A closed symplectic four-manifold  has
 almost K\"ahler metrics of negative scalar curvature}
%\email{jskim@ccs.sogang.ac.kr}

\author{
%Yutae Kang\\ E-mail: lubo@sogang.ac.kr \and
 \ Jongsu Kim
\thanks{Supported by grant No. R14-2002-044-01000-0(2002) from
the Korea Science and Engineering Foundation. Keywords: Almost
K\"ahler metric, symplectic structure, scalar curvature.
MS Classification(2000): 53C15, 53C20, 53D35} \\
%E-mail: jskim@sogang.ac.kr \\
\and
Department of Mathematics, Sogang University\\
%Sinsu-dong 1, Mapo-gu  \\
Seoul, 121-742, KOREA }
%\date{May. 3, 2005}

\maketitle

\begin{abstract}
We show that every closed symplectic four-dimensional manifold $(M,
\omega)$ admits an almost \k metric of negative scalar curvature
compatible with $\omega$.
\end{abstract}

\section{Introduction}

In Riemannian geometry an interesting problem related to a curvature
is whether every smooth manifold of dimension bigger than 2 admits a
metric with the curvature negative. Lohkamp's beautiful argument in
\cite{Lo1} resolved this for the Ricci curvature case, while the
relatively simpler case of scalar curvature had been understood
earlier, see \cite[Chap. 4]{Be}.

%An almost K\"{a}hler structure on a smooth manifold $M$ is a triple $(g,\omega,J)$ where $\omega$ is a symplectic form, $J$ an
%almost complex structure, and $g$ is an $\omega$-compatible Riemannian metric, i.e. $\omega(x,y) = g(Jx,y)$ for tangent
%vectors $x,y$ to $M$. A manifold with an almost \k  structure is called an almost \k manifold.

\medskip
Meanwhile, the extensive study on symplectic manifolds in the last
decades has broadened the knowledge on the field and supplied large
families of closed symplectic manifolds \cite{GS, Gr}. Therefore the
study of Riemannian metrics compatible to a symplectic structure,
so-called almost \k metrics, re-emerges as an interesting and
promising subject.  Remarkably, the study of almost \k metrics in
Riemannian geometry has a long history, from the works of 1960's or
earlier \cite{Y}
 to recent works in several directions \cite{ADK}, \cite{I}, \cite{KS}, etc..

\medskip
Motivated by this, in this paper we study the existence of negative
scalar-curved  almost \k metrics on a symplectic manifold.
In the
 general Riemannian  metric case
afore-mentioned, the main argument dealing with scalar curvature was
to use a conformal deformation. For almost \k metrics the latter can
not be used because it does not preserve the symplectic structure.
So even the scalar curvature question does not seem trivial in
almost \k case.

Our argument is based on the Lohkamp's framework in \cite{Lo1} and in its course there were two
main ingredients to  proving the following main theorem in this
paper;

\begin{thm} \label{th1}
Every closed symplectic four-dimensional manifold $(M, \omega)$
admits an $\omega$-compatible almost \k metric of negative scalar
curvature.
\end{thm}

The first ingredient is to find an almost \k {\it island} metric,
i.e. an almost \k metric on ${\Bbb R}^n$  which has the scalar
curvature negative on a pre-compact open subset and is Euclidean
outside of its closure.
%This part is done in \cite{KK}.

 The second one is to find an almost \k deformation to
decrease the scalar curvature outside copies of the island metrics
which are embedded into the (enlarged) almost \k manifold by a
metric surgery. This deformation corresponds to the conformal
deformation \cite[p. 669]{Lo1} of Lohkamp's.

As our almost \k deformation completely avoids any conformal
deformation, it brings about new technical problems.
% and significant differences from.
One remarkable difference from Lohkamp's argument is that all the
metrics involved, including the island metrics, should be locally
close to Euclidean metrics. It is because we resorted to specific
local orthonormal frame fields $\{ \omega_i \}$ in defining the
deformation, and these fields are manageable if the metrics are
close to Euclidean ones. This technical aspect is the reason why we
restricted to real 4-dimension in this article, though we believe
the higher dimensional case should be provable by some modification.

\medskip
One may view the sort of result like Theorem \ref{th1} as indicating
that the space of almost \k metrics on a symplectic manifold may be
big enough to contain rich geometric features. And one may speculate
on any implication of almost \k geometry to the topology of
symplectic manifolds.

\par

\bigskip
The paper is organized as follows. In section 2 we provide
preliminaries for almost \k metrics, an estimate for injectivity
radius and the Besicovitch covering. In section 3 we describe the
local orthonormal frame fields $\{ \omega_i \}$ of almost \k metrics
which are close to Euclidean metrics. In section 4 we explain the
almost \k deformation depending on $\{ \omega_i \}$ which decreases
the scalar curvature away from a compact set. In section 5 we
describe an {\it island} metric and modify it with the result of
section 4. In section 6 we do a surgery of embedding many copies of
metrics of section 5 into the given almost \k manifold. In section 7
we apply the deformation of the section 4 to the resulting manifold
of section 6 and prove the main theorem. There is Appendix which
proves some estimates on  $\{ \omega_i \}$ and derives a scalar
curvature formula.

\medskip
{\bf Acknowledgement.} The author would like to thank  Yutae Kang
for many helpful discussions on this work and thank the referee for
his suggestions which made this paper more readable.

\section{Preliminaries}

In this section we explain definitions and formulas which will be needed in later sections.

\indent An {\it almost-K\"{a}hler} metric  is a Riemannian metric
$g$ compatible with a symplectic structure $\omega$ on a smooth
manifold, i.e. $ \omega(x, y)= g(Jx, y) $ for an almost complex
structure  $J$, where $x, y$ are tangent vectors at a point of the
manifold. We often denote it by the triple  $(g, \omega, J)$ or the
couple $(g, \omega)$ to specify what  $\omega$ or $J$ is. Notice
that any one of the pairs $(g, \omega)$, $(\omega, J)$, or $(g, J)$
determines the other two. For a fixed $\omega$, we shall call a
metric $g$ $\omega$-almost \k if $g$ is compatible with $\omega$ and
denote by $\Omega_{\omega}$
 the set of all smooth $\omega$-almost \k metrics.
An almost-\k metric $(g, \omega, J)$ is \k if and only if  $J$ is
integrable.

An almost complex structure $J$ gives rise to a type decomposition
of any symmetric (2,0)-tensor $h$ which  can be written as $h = h^+
+ h^-$, where $h^+$ and $h^-$ are symmetric 2-tensors defined by
$h^+ (x, y) = \frac{1}{2} \{ h(x, y) + h(Jx, Jy) \}$ and  $h^- (x,
y) = \frac{1}{2} \{ h(x, y) - h(Jx, Jy) \}$ for two vectors $x$ and
$y$ tangent to $M$. A symmetric $(2,0)$-tensor $h$ is called  {\it
$J$-invariant} or {\it $J$-anti-invariant} if $h = h^+$ or $h = h^-$
respectively.

The set $\Omega_{\omega}$ is naturally an infinite dimensional
Fr\'echet manifold. %According to Blair \cite{Bl},
For a smooth curve $g_t$ in $\Omega_{\omega}$ with the corresponding
curve $J_t$ of almost complex structures, $h = \frac{d g_t}{dt}
|_{t=0}$ is $J_0$-anti-invariant. Conversely, for $g$ in
$\Omega_\omega$ with the corresponding $J$, any $J$-anti-invariant
symmetric $(2,0)$-tensor field $h$  is tangent to a smooth curve in
$\Omega_{\omega}$. More precisely, $t \mapsto g e^{ht}$ is such a
smooth curve in
$\Omega_\omega$,  %with corresponding $J \cdot e^{ht}$
where  $g e^{ht}$ is defined by $g e^{ht} (x, y) = g(x, e^{(g^{-1}h)t} y)$,
i.e.
\begin{equation} \label{a001}
g e^{ht} (x, y)= g(x,y) + g(x, \sum_{k=1}^{\infty} {t^k (g^{-1}{h})^k
(y) \over k!}),
\end{equation}
\noindent
 where $g^{-1} {h}$  is the lifted $(1,1)$-tensor of $h$
with respect to $g$. Given $g \in \Omega_{\omega}$ with
corresponding $J$, any other metric $\tilde{g}$ in $
\Omega_{\omega}$ can be expressed as

\begin{equation}
\tilde{g}= g e^h , \label{geh}
\end{equation}
where $h$ is a $J$-anti-invariant symmetric $(2,0)$-tensor field
uniquely determined. See \cite{Bl} for details.

\bigskip
\noindent
 We recall local lower bound estimates on the injectivity
radius of a manifold. For  $p $ in a complete Riemannian manifold
$(N, g)$, let $l_p$ denote the length of the shortest geodesic loop
based at $p$. If the sectional curvature $K_N \leq \Lambda^2$ on
$B^g_r (p)$, the open geodesic ball of radius $r$ centered at $p$
w.r.t. $g$, and $r \leq \frac{1}{ 2} l_p$, then the injectivity
radius $\rm{inj}_g (p)$ of $(N, g)$ at $p$ satisfies $\rm{inj}_g (p)
\geq \rm{min} (r, \frac{\pi}{ \Lambda})$ by an estimate of
Klingenberg \cite[Lemma 1]{Kl} and the fact that if a point $q$ is
conjugate to $p$ then their distance $d(p, q) \geq
\frac{\pi}{\Lambda}$ ,  \cite[p. 30]{CE}.

\smallskip \noindent
Now suppose $K_N \leq \Lambda^2$ on $B^g_r (p)$,  where $r \leq
\frac{\pi} {\Lambda}$. Choose any numbers $r_0$ and $s$ such that
$r_0 +2s \leq r$, $r_0 \leq \frac{r }{4}$. Then one has the
estimate due to Cheeger-Gromov-Taylor \cite[Theorem 4.3]{CGT}

\begin{equation}
l_p \geq r_0 [1 + \frac{v^{- \Lambda}_{r_0 + s} } {v^g_s (p)} ]
^{-1}. \label{int}
\end{equation}
 \noindent where $v^g_s (p)$ is the $g$-volume of $B^g_s
(p)$ and $v^{- \Lambda}_{r_0 + s} $ is the volume of the geodesic
ball of radius $r_0 + s$ in the space form of constant curvature $-
\Lambda$.

\bigskip
Next we recall the Besicovitch covering \cite[Prop 4.4]{Lo1};
\begin{lemma} \label{7}
 Let $(M, h)$ be a Riemannian manifold with the
injectivity radius $\rm{inj}_h \geq 100$ and the sectional curvature
$|K_h| < 1$. Then there are constants $\kappa \in Z^{>0}$ and $m_0
\geq  1$ depending only on
 the dimension of $M$ such that

\rm{(i)} For each  $m \geq m_0$ there is a set $A=A(m) \subset M$
with $d_{m^2 h} (a, b) > 5$ for $a \neq b, a, b \in A(m)$.

\smallskip
\rm{(ii)} $ {\cal A}(m) =  \{ \overline{B^{m^2 h}_5 (a)}  | \mbox{
}a \in A(m) \} $
 is a (closed) covering of $M$ and splits into $\kappa$ disjoint
 families ${\cal B}_j$ with
 \( \overline{B^{m^2 h}_{10} (a)}  \cap  \overline{B^{m^2 h}_{10} (b)}
  = \phi    \)
 if
$\overline{B^{m^2 h}_5 (a)}$  and $\overline{B^{m^2 h}_5 (b)}$
belong to the same ${\cal B}_j$.
\end{lemma}

\bigskip
We shall use the following notation. When $h$ is a tensor field on a
domain $V$, we denote the $C^k$ norm of
 $h$ with respect to a metric $g$ on $V$  by $\| h
 \|_{{C^k_{g}}(V)}$. Or sometimes we express it by writing $``\| h
 \|_{{C^k_{g}}}$ on $V"$. Here $\| h
 \|_{{C^k_{g}}(V)}= \sum_{i=0}^k \| (\nabla_g)^i h \|_{{C^0_{g}}(V)}
 $ and $\nabla_g$ is the Levi-Civita connection of $g$.

\section{Almost-K\"{a}hler metrics which are close to Euclidean metrics}

\indent We start with a Euclidean metric $g_0$ on ${\Bbb R}^4$ or an
open subset ${\Omega}$. We assume that there exists a point $p$ in $
\Omega$ and
 that $B^{g_0}_{10^4}(p)$, the open $g_0$-ball
centered at $p$ with radius $10^4$, is in $\Omega$. Now suppose that
there is another Riemannian metric $g$ on $ \Omega$ and $\|g- g_{0}
\|_{{C^0_{g_{0}}}( B^{g_0}_{10^4}(p))} \leq \varepsilon $, where
$\varepsilon$ is a small positive number, say $\varepsilon <
10^{-4}$.

 For another point $q \in B^{g_0}_{10^4} (p)$ and a minimal
geodesic $\gamma_{g_{0}}(t)$ of $g_{0}$ connecting  $\gamma_{g_{0}}
(0)= p$, $\gamma_{g_{0}} (t_0)= q $, we have
$$ d_{g_{0}} (p, q) = \int_0^{t_0}  |\gamma_{g_{0}}^{\prime}(s)|_{{g_{0}}} ds \geq
    \int_0^{t_0} \frac{|\gamma_{g_{0}}^{\prime}(s)|_{g}}{\sqrt{1+ \varepsilon}}   ds
    \geq  \frac{d_{g} (p, q) }{\sqrt{1+ \varepsilon}}  ,$$
\noindent where $ d_g (p, \cdot)$ and $ d_{g_{0}} (p, \cdot)$ are
the distance functions from $p$ induced by $g$ and $g_{0}$,
respectively.
 \noindent To prove the inequality in the opposite
direction, we first show $B^{g}_{9999} (p) \subset B^{g_0}_{10^4}
(p)$. Let's consider a $g$-geodesic $\gamma(t)$ with $\gamma(0)=p$
and $|\gamma^{\prime}(t)|_g \equiv 1$. For a small number $\delta$
with $0 < \delta \ll \varepsilon$, let $t_1 >0$ be the first
positive value of $t$ such that $\gamma(t)$ is on the boundary of $
B^{g_0}_{10^4 - \delta} (p)$.
 Then $$ 10^4 - \delta = d_{g_0}(p, \gamma(t_1)) \leq \int_0^{t_1} |
\gamma^{\prime}{(s)} |_{g_0} ds \leq \frac{1}{\sqrt{1-\varepsilon}}
\int_0^{t_1} | \gamma^{\prime}{(s)} |_{g} ds =
\frac{1}{\sqrt{1-\varepsilon}}  t_1. $$

As $9999 < (10^4 -\delta ){\sqrt{1-\varepsilon}} \leq t_1$, this
implies $B^{g}_{9999} (p) \subset B^{g_0}_{10^4} (p)$. Now for $q
\in B^{g}_{9999} (p) $ and a minimal geodesic $\gamma_g (t)$ of $g$
connecting  $\gamma_g (0)= p$, $\gamma_g (t_2)= q $, we have
$$ d_g (p, q) = \int_0^{t_2} |\gamma_g^{\prime}(s)|_{g} ds \geq
\int_0^{t_2} \sqrt{1- \varepsilon} \cdot
|\gamma_{g}^{\prime}(s)|_{g_0} ds  \geq \sqrt{1- \varepsilon} \cdot
d_{g_{0}} (p, q) .$$

\noindent In sum we proved;
 \begin{lemma} \label{lam1}
Suppose that a Euclidean metric $g_0$ and a Riemannian metric $g$ on
 $B^{g_0}_{10^4}(p) \subset {\Bbb R}^4$ satisfies
  $\|g- g_{0}
\|_{{C^0_{g_{0}}}( B^{g_0}_{10^4}(p))} \leq \varepsilon < 10^{-4}$.
Then, $B^{g}_{9999} (p) \subset B^{g_0}_{10^4} (p)$ and
 for any $ q \in B^{g}_{9999} (p)$,
we have
\begin{equation} \label{dist}
\sqrt{1- \varepsilon} \cdot d_{g_{0}} (p, q)  \leq    d_g (p, q)
  \leq  \sqrt{1+ \varepsilon}  \cdot d_{g_{0}} (p,
q).
\end{equation}
\end{lemma}

One can prove the following lemma from (\ref{int}) and (\ref{dist}).
We omit the detail of the proof which is straightforward.
%One just note that the estimates of Klingenberg's and Cheeger-Gromov-Taylor's are from purely
%local arguments so that they can be applied for the metrics of the next lemma.
 \begin{lemma} \label{lamb}
There exists $\lambda_1 < 10^{-4}$ such that if $\|g- g_{0}
\|_{{C^2_{g_{0}}}(B^{g_{0}}_{10^4}(p))} \leq \lambda_1$, then $g$
has the injectivity radius at $p$,
 $\rm{inj}_g (p) \geq 200$.
\end{lemma}

\bigskip
Consider the symplectic structure $\omega= dx_1 \wedge dx_{2} + dx_3
\wedge dx_{4}$ on ${\Bbb R}^{4}= \{(x_1, x_2, x_3, x_4)  | \mbox{ } x_i \in  {\Bbb R} \}$
together with the compatible Euclidean metric $g_0 = \sum_{i=1}^{4}
dx_i \otimes dx_i$ and the (almost) complex structure $J_0$. We
shall deform an $\omega$-almost \k metric $g$ on  $\Omega$
containing $ B^{g_0}_{10^4} (p)$, satisfying $\| g- g_{0}
\|_{{C^l_{g_{0}}}(B^{g_{0}}_{10^4}(p))} \leq \varepsilon < 10^{-4}$,
where $ l \geq 4$ and $\varepsilon$ is a positive number to be determined later. The
Euclidean metric $({\Bbb R}^4, g_0, \omega, J_0) $ has an
orthonormal co-frame field $\omega_1^0= dr_{g_0}$, $\omega_2^0 = J_0
(dr_{g_{0}}) = r_{g_0} \sigma_1$, $\omega_3^0 = r_{g_0} \sigma_2$,
 $\omega_4^0 =  r_{g_0} \sigma_3 = J_0 (\omega_3^0)$, where $r_{g_0}(\cdot)
 = d_{g_0}(p, \cdot)$ and  $\sigma_1,  \sigma_2, \sigma_3$ is the canonical
 orthonormal co-frame field  of the standard metric on $S^3$ which satisfy $ d \sigma_1= 2\sigma_2 \wedge \sigma_3,
  d \sigma_2=2 \sigma_3 \wedge \sigma_1$ and $d \sigma_3=2 \sigma_1 \wedge \sigma_2$.
 We shall find an
orthonormal co-frame field of $g$ near $p$:
%on $B^{g_{0}}_{11}(0)$:
$\omega_i$, $i=1, 2, 3, 4$ which are close to $\omega_i^0$. On a
ball $ B^g_{100}(p) -\{p\} \subset B^{g_{0}}_{101}(p)$ we define formally
\begin{align}
\omega_1 &= dr_g, \hspace{.5in } \omega_2 = J_g (dr_g), \nonumber \\
  \omega_3 &=\frac{
\omega_3^0 -g(\omega_3^0, \omega_1)\omega_1 - g(\omega_3^0,
\omega_2)\omega_2}{|    \omega_3^0 -g(\omega_3^0, \omega_1)\omega_1
- g(\omega_3^0, \omega_2)\omega_2 |_g }, \hspace{.2in } \omega_4 =
J_g (\omega_3). \label{oms}
\end{align}

\noindent where $r_{g}(\cdot)
 = d_{g}(p, \cdot)$ and $J_g$ is the almost complex structure determined by
$\omega$ and $g$. We need the following proposition on some
properties of $\{ \omega_i \} $ whose proof is postponed to the
Appendix;

\begin{prop} \label{lemom}
There exist $\lambda_2$ with $0 < \lambda_2 \leq \lambda_1$ and a
positive constant $C$ such that  if $\|g- g_{0} \|_{{C^l_{g_{0}}}}
\leq \lambda  \leq \lambda_2 $ on $B^{g_{0}}_{10^4}(p) \subset {\Bbb
R}^4$, where $l$ is an integer $4 \leq l \leq 2 \kappa +4$ and
$\kappa$ is the constant (in dimension four) of the Besicovitch
covering Lemma \ref{7}, then
 $\omega_1, \omega_2, \omega_3, \omega_4$ in the formula (\ref{oms}) are well defined and satisfy on $B^g_{100}(p) - \{
 p
 \}$;
\begin{align}\label{prop1}
 \| (\nabla^{g_0})^j \omega_i - (\nabla^{g_0})^j \omega_i^0 \|_{{C^0_{g_{0}}}}  <
 \frac{C \cdot {\lambda}}{r^j_{g_{0}}},  \hspace{0.2in} \mbox{  for }
 j =0, 1, \cdots, l-2.
\end{align}
\noindent
 where $  \nabla^{g_0}$ is the Levi-Civita connection of $g_0$.
\end{prop}

\noindent For the above metric $g= \sum_{i=1}^{4} \omega_i \otimes
\omega_i$ on $ B^{g}_{100}(p)$, consider the connection 1-forms
${\omega}_{ij}$ w.r.t. the co-frame $\omega_i$: $d {\omega}_{i } =
 \sum_{j=1}^{4}{\omega}_{ij} \wedge {
\omega}_{j}$, with $ {\omega}_{ij}= -{\omega}_{ji}$.
%Of course, these equations determine ${\omega}_{ij}$ uniquely. Indeed,
Writing ${\omega}_{ij}$ as $ {\omega}_{ij}=\sum_{k=1}^{4} a_{ijk}
\cdot {\omega}_{k}$, we get $d {\omega}_{i } = \sum_{j < k}
(a_{ikj}- a_{ijk}) {\omega}_{j} \wedge {\omega}_{k}$. So we get
$2a_{ijk} = \langle d {\omega}_{k}, {\omega}_{i} \wedge {\omega}_{j}
\rangle_g -\langle d {\omega}_{i}, {\omega}_{j} \wedge {\omega}_{k}
\rangle_g -\langle d {\omega}_{j}, {\omega}_{k} \wedge {\omega}_{i}
\rangle_g$, where $\langle \cdot ,\cdot \rangle_g$ is the naturally
induced metric on the space of 2-forms from $g$.
%Now, three values $a_{ijk}$, $a_{jki}$ and $a_{kij}$ can be determined at once
%as we know $a_{ijk}- a_{ikj} =a_{ijk}+ a_{kij}$, $a_{jik}- a_{jki}=
%-a_{ijk}- a_{jki}$ and $a_{kij}- a_{kji}=a_{kij}+ a_{jki} $.

%with respect to the basis ${\omega}_{1}$,${\omega}_{2}$, $\omega_{3}$, $\omega_{4}$.
For the Euclidean  metric $({\Bbb R}^4, g_0, \omega, J_0) $ with the
above co-frame $\{ \omega^0_i \}$,
%$\omega^0_1= dr$, $\omega^0_2 = J_0 (dr) = r \sigma_2$, $\omega^0_3 = r \sigma_3$, $\omega^0_4 =  r \sigma_4$,
the $a_{ijk}$'s are all zero except

\begin{equation} \left\{ \begin{array}{ll}
a_{122} &= a_{133}= a_{144}= -a_{212} =- a_{313}= -a_{414}= \frac{1}{r_{g_0}},  \\
a_{234} &= a_{342}= a_{423}= -a_{324} =- a_{432}= -a_{243}=
-\frac{1}{r_{g_0}}.
          \end{array}
          \right. \label{aijk} \end{equation}
For later we need the following estimate.
%We write its proof in the Appendix.

\begin{lemma} \label{lemaijk}
There exists $\lambda_3$ with  $ \lambda_3 \leq \lambda_2 $ such
that
 if $\|g- g_{0} \|_{{C^4_{g_{0}}}} \leq
\lambda_3$ on $B^{g_{0}}_{10^4}(p) \subset {\Bbb R}^4$, then  for
$i,j,k, l=1, 2,3,4$, the  $a_{ijk}$ of the $\omega_i$'s in the
formula {\rm (\ref{oms})} satisfy:
$$ |a_{ijk}| < \frac{2}{r_{g}}, \mbox{   } |a_{ijk, l}| <
\frac{2}{r_{g}^2} \mbox{ on } B^g_{100}(p) - \{p \},$$
\end{lemma}
\noindent where $a_{ijk,l}$  is defined in the formula $d(a_{ijk}) =
\sum_{l=1}^{4} a_{ijk,l} \omega_l$.

\bigskip \noindent {\bf Proof:}
%of Lemma \ref{lemom}
%Let $e_i$, $i=1,2,3,4,$ be the dual frame-field of $\{ \omega_i \}$. $a_{ijk}$ is defined in the system
%of equations $ d \omega_i = \sum_{j} \omega_{ij} \wedge \omega_j = \sum_{j<k} (a_{ijk} - a_{ikj}) \omega_k
%\wedge \omega_j$ so that $a_{ijk}$ can be expressed as a linear combination of  $d \omega_{\alpha}
%( e_{\beta}, e_{\gamma})$. For instance, $a_{123}$ is a linear combination of $d
%\omega_1( e_2, e_3)$, $d \omega_2( e_3, e_1)$ and $d \omega_3( e_1, e_2)$.
One needs a straightforward computation:  show
 $|\langle d{\omega}_{i}, {\omega}_{j} \wedge {\omega}_{k}
\rangle_g - \langle d {\omega}_{i}^0, {\omega}_{j}^0 \wedge
{\omega}_{k}^0 \rangle_g| \leq \frac{{\rm constant} \cdot
\lambda}{r_{g_0}}$  from (\ref{prop1}). Then $a_{ijk} =
\frac{1}{2}(\langle d {\omega}_{k}, {\omega}_{i} \wedge {\omega}_{j}
\rangle_g -\langle d {\omega}_{i}, {\omega}_{j} \wedge {\omega}_{k}
\rangle_g -\langle d {\omega}_{j}, {\omega}_{k} \wedge {\omega}_{i}
\rangle_g)$ satisfies $|a_{ijk} -a_{ijk}^0| \leq \frac{{\rm
constant} \cdot \lambda}{r_{g_0}}$. From (\ref{aijk}), $ |a_{ijk}| <
\frac{2}{r_{g}}$ follows. One can do similarly for $|a_{ijk, l}|$.

\section{Scalar curvature diffusion}

We shall perturb locally an almost \k metric $g$ which satisfies the
hypothesis of Lemma \ref{lemaijk} and is written as $g= \sum_{i=1}^4
\omega_i \otimes \omega_i$ on $ B^g_{100}(p) - \{p \}$, where
$\omega_i$'s are as defined in (\ref{oms}).

 We use the functions of \cite{Lo1};
 $ f_{d,s}(t) \in C^{\infty}(\Bbb{R}, \Bbb{R}^{\geq 0} )$ for $d, s > 0$ by
 $f_{d,s}  = s \cdot \rm{exp}(- \frac{d}{t})$ on $\Bbb{R}^{>0}$ and
 $f_{d,s}  =  0$ on  $\Bbb{R}^{\leq 0}$.
%and let $f_{d,s}\in C^{\infty}(\Bbb{R}^{2n}, \Bbb{R}^{\geq 0} ) $ be $f_{d,s}(t,\cdots) =  F_{d,s}(t)$.
Also choose an $h \in C^{\infty}(\Bbb{R}, [0,1])$ with $ h= 0$ on
$\Bbb{R}^{\geq 1}$, $h=1$ on $\Bbb{R}^{\leq 0}$ and
$h_{\epsilon}^{b}(t) = h( \frac{1}{\epsilon}(t-b)), $ $ b > 0$, $
\epsilon >0$.

\noindent
 Now we define a perturbed metric for $b+ \epsilon <c \leq 9$ ;
% \frac{1}{\alpha^2}-1 = Y(r) = h_{\epsilon}^{b} \cdot f_{d,s} and so from (\ref{gtilde})
 \begin{align}
 g_{d,s}^{b, \epsilon,c}  =& \frac{ \omega_1 \otimes
\omega_1}{\{ 1+ h_{\epsilon}^{b}(c-r_g) \cdot f_{d,s}(c-r_g) \}} +
\{ 1+ h_{\epsilon}^{b}(c-r_g) \cdot f_{d,s}(c-r_g)\}
  \omega_2 \otimes \omega_2 \nonumber  \\
  &+ \omega_3 \otimes \omega_3
 +\omega_4 \otimes \omega_4. \label{gbe}
 \end{align}

\noindent Then $ g_{d,s}^{b, \epsilon, c}$ is a smooth
$\omega$-almost \k metric on $B^{g_0}_{10^4} (p)$. We get

\begin{lemma} \label{5}
Suppose that an
$\omega$-almost \k metric $g$ on a domain $\Omega \subset
\Bbb{R}^{4}$ satisfies $\|g- g_{0}\|_{{C^l_{g_{0}}}
(B^{g_0}_{10^4}(p))} \leq \lambda_3$ where $ B^{g_0}_{10^4}(p)
\subset \Omega$ and $l$ is an integer $4 \leq l \leq 2 \kappa +4$. For each $a, b, \epsilon, c
>0$ such that $ 9 \geq c
> b+\epsilon >b >a >0$, consider the orthonormal co-frame $\omega_i$ as defined in \rm{(\ref{oms})}
and define $g_{d,s}^{b, \epsilon, c}$ on $ B^{g_0}_{10^4} (p)$ as in
(\ref{gbe}). Then the following statements hold;

 \noindent
(i) For each $D >0$, there is a constant $\alpha = \alpha(b,
\epsilon, c, D)$ such that for $d \geq D$ and $s \in (0,1]$,
$\|g_{d,s}^{b, \epsilon, c}-g\|_{{C^{l-2}_{g_{0}}} ( B^{g_0}_{10^4}
(p))} < s \cdot \alpha$,

\noindent (ii) For each $\varepsilon >0$ there exists a $d= d
({\varepsilon}, c)$ such that for every $d \geq d({\varepsilon},c)$
and each $s \in (0,1]$, $\|g_{d,s}^{b, \epsilon, c}-g
\|_{{C^{l-2}_{g_{0}}} (B^{g_0}_{10^4} (p))} < \varepsilon$ holds.
\end{lemma}
\begin{proof}
  We apply the estimates (\ref{prop1}) to
  %for the derivatives of $r_g$  (for $\omega_1 = dr_g$) in  Proposition \ref{lemom} to
  $$g_{d,s}^{b,
\epsilon, c}-g = - \frac{h_{\epsilon}^{b}(c-r_g) \cdot
f_{d,s}(c-r_g)}{1 + h_{\epsilon}^{b}(c-r_g) \cdot f_{d,s}} \omega_1
\otimes  \omega_1 +
 h_{\epsilon}^{b}(c-r_g) \cdot f_{d,s}(c-r_g)\omega_2 \otimes  \omega_2$$
 which has support in $ c-b- \epsilon \leq r_g \leq c$.
The computation is elementary and we omit it. Here one may use $ \|
(\nabla^{g_0})^j dr_{g_0} \|_{{C^0_{g_{0}}}}  <
 \frac{{\rm constant}}{r^j_{g_{0}}}$, for $0 \leq j \leq l-2$;  see the
formula (\ref{35}) in the Appendix.
\end{proof}

Note that $\lambda_3$ is a constant already determined, so in the
statement of Lemma \ref{5} we do not specify the dependence of
$\alpha$ or $d$ on $\lambda_3$.

For the metric $g_{d,s}^{b, \epsilon, c}$, we compute its scalar
curvature using the orthonormal co-frame field
$\tilde{{\omega}_{1}}= {{\omega}_{1} \over {\sqrt{ 1+
h_{\epsilon}^{b} \cdot f_{d,s}(c-r_g) }}} $,
$\tilde{{\omega}_{2}}=\sqrt{ 1+ h_{\epsilon}^{b} \cdot
f_{d,s}(c-r_g) } \cdot {\omega}_{2}$,
$\tilde{{\omega}_{3}}={\omega}_{3}$ and
$\tilde{{\omega}_{4}}={\omega}_{4}$ and the formula $d
\tilde{{\omega}}_{i } =
 \sum_{ j=1}^{4} \tilde{{\omega}}_{ij} \wedge \tilde{{\omega}}_{j}$.

\noindent Setting  $Y=Y(r_g) = h_{\epsilon}^{b}(c-r_g) \cdot f_{d,s}
(c-r_g)$,
%$Y^{'}=  Y^{'}(r_g)$   and $Y^{''}=Y^{''}(r_g) $,
the scalar curvature of $g_{d,s}^{b, \epsilon, c}$ is computed as
follows, see Appendix for details:
\begin{align}
{ s(g_{d,s}^{b, \epsilon, c}) }
%=& 2(R_{2112} + R_{3113} + R_{4114} + R_{3223} + R_{4224} + R_{4334}) \\
          =&-Y^{''}- Y^{'} \cdot (3a_{122} +
2 a_{133} + 2 a_{144})-2(Y+1)\{\sum_{i=2}^{4} (a_{1ii,1}
           +a_{1ii}^{2}) \nonumber  \\
          & +a_{134}^2+a_{122}a_{133}+a_{122}a_{144}+ a_{133}a_{144}
+{1 \over 4}(a_{243}-a_{234})^{2}\} \nonumber \\
          & -{2 \over {Y+1}}\{a_{233,2}+ a_{244,2}+a_{233}^{2}+a_{244}^{2}+a_{233}a_{244}
          +{1 \over 4}(a_{243}+a_{234})^{2} \} \nonumber \\
          & - {1 \over 2}\sum_{i=3}^4\{Ya_{2i1}-(Y+2)a_{21i}\}^{2}  + I_0,
          \label{sgt}
\end{align}
\begin{align}
\mbox{where } \frac{I_0}{2}= &
a_{232,3}+a_{242,4}+a_{343,4}+a_{434,3}-
          a_{232}^{2}-a_{242}^{2}-a_{343}^{2}-a_{434}^{2}+a_{242}a_{433} \nonumber\\
          &+a_{232}a_{344}+a_{342}a_{234}-a_{342}a_{243}+{1 \over 2}a_{234}^{2}+{1 \over 2}a_{243}^{2}
          -a_{243}a_{234}. \nonumber
\end{align}
\noindent The scalar curvature of $g$ can be obtained by putting
$Y=0$ into (\ref{sgt}). We are concerned with their difference:

\begin{equation}
{ s(g_{d,s}^{b, \epsilon, c})-s(g) } = - Y^{''}-  Y^{'}(3a_{122} + 2
a_{133} + 2 a_{144})- Y \cdot I_1 \label{sdif}
\end{equation}
\noindent where
 \begin{align*}
 {I_1 \over 2}  = &  {1 \over 4}(a_{243}-a_{234})^{2} +\sum_{i=2}^{4} (a_{1ii,1} +a_{1ii}^{2})
 +a_{134}^2+a_{122}a_{133}+a_{122}a_{144} + a_{133}a_{144}\\
 + &{1 \over 4} \{Y(a_{231}^2 + a_{241}^2)-2(Y+2)(a_{213} a_{231} +
a_{214} a_{241}) +(Y+4)(a_{213}^2  +  a_{214}^2 )\}   \\
 -& ({1 \over {Y+1}})\{a_{233,2}+
a_{244,2}+a_{233}^{2}+a_{244}^{2}+a_{233}a_{244}+{1 \over
4}(a_{243}+a_{234})^{2} \}.
  \end{align*}

Due to the estimate of Lemma \ref{lemaijk}, we can prove that the
above almost \k deformation diffuses the scalar curvature,
analogously to \cite[Proposition 2.1]{Lo1}

\begin{lemma} \label{6}
\smallskip

 \noindent
Suppose that an almost \k metric  $g$ on a domain $\Omega$ in  $
\Bbb{R}^{4}$ satisfies $\|g- g_{0}\|_{{C^4_{g_{0}}} (B^{g_0}_{10^4
}(p))} \leq \lambda_3
 < 1$. For each $a, b,  \epsilon, c >0$ such that $ 9 \geq c > b+\epsilon >b >a >0$,
consider the orthonormal co-frame $\omega_i$ as in (\ref{oms}) and
define $g_{d,s}^{b, \epsilon, c}$ on $ B^{g_0}_{10^4} (p)$ as in
(\ref{gbe}). There are constants $\gamma = \gamma(a, b, c)>0$
%which depend only on $a,b$
 such that for $d \geq \gamma$, $s \in [0,1]$ the following statements hold:

\rm{(i)} $ g_{d,s}^{b, \epsilon, c} \equiv g $ on $B^g_{c-b -
\epsilon}(p)$,
%$ \Bbb{R} \setminus (0, b+ \epsilon) \times \Bbb{R}^n $.

\rm{(ii)} $s({g_{d,s}^{b, \epsilon, c}})-s(g) \leq 0$ on
$B^{g_0}_{10^4} (p) \setminus  B^g_{c-b}(p) $,
%$(0,b] \times \Bbb{R}^n$,

\rm{(iii)} $s({g_{d,s}^{b, \epsilon, c}})-s(g) \leq -s \cdot e^{-
\frac{d}{a}}$ on   $   B^g_{c-a}(p)  \setminus B^g_{c-b}(p)  $,
%$(a,b] \times \Bbb{R}^n$.

\end{lemma}

\begin{proof}
(i) is clear from the definition. To show (ii) and (iii), one needs
to consider terms in the right hand side of (\ref{sdif}). From the
assumption on $g$ and Lemma \ref{lemaijk} we have $|a_{ijk}| <
\frac{2}{r_g}$ and $|a_{ijk, l}| < \frac{2}{r_{g}^2}$. Note that  $
0 \leq Y(r_g) \leq 1$ for $0 \leq s \leq 1$ and on
 $B^{g_0}_{10^4} (p)$, for any $d >0$.
 So we have $| I_1| \leq
\frac{\nu}{r_g^2}$ for a constant $\nu$. By some computation
(\cite[Lemma 1.2 (i)]{Lo1}),  there is $d_0 (b)>0 $ such that for $
d \geq d_0(b)$,
 $f^{(k)}_{d,s} > 0$ on $(0, b]$ for $k=0,1,2,3$. We assert
that there is $d_0 (a,b,c) > d_0 (b)$ such that for $d \geq d_0
(a,b,c)$ and every $s >0$ the following two inequalities hold:
$$ f^{''}_{d,s} -  f^{'}_{d,s} \frac{c_1 }{c-t} -  f_{d,s} \cdot \frac{\nu }{(c-t)^2} \geq 0,  \hspace{0.1in} \mbox{ on } (0, b],
$$
$$  \hspace{0.5in} f^{''}_{d,s} -  f^{'}_{d,s} \frac{c_1}{c-t} -  f_{d,s} \cdot \frac{\nu }{(c-t)^2} \geq s \cdot e^{-
\frac{d}{a}}, \hspace{0.1in} \mbox{  on }  (a, b], \hspace{0.2in}$$
where  $c_1$ is a constant such that $|3a_{122} + 2 a_{133} + 2
a_{144}| \leq \frac{c_1}{r_g}$. The proof of this assertion is
similar to that of Lemma 1.2 (ii) in \cite{Lo1}, so we omit it.

\noindent Then we get on  $B^g_{c} (p) \setminus B^g_{c-b}(p)$,
where $ h_{\epsilon}^{b}(c-r_g) \equiv 1$,
\begin{align}
 Y^{''} + Y^{'}(3a_{122} + & 2 a_{133} + 2 a_{144})+ Y \cdot
I_1 \nonumber \\
\geq & f^{''}_{d,s}(c-r_g) - f^{'}_{d,s}(c-r_g) \frac{c_1 }{r_g} -
f_{d,s}(c-r_g) \cdot \frac{\nu }{(r_g)^2}  \nonumber \\
= & f^{''}_{d,s}(t) -  f^{'}_{d,s}(t) \frac{c_1 }{c-t} -  f_{d,s}(t)
\cdot \frac{\nu }{(c-t)^2},  \nonumber
\end{align}
where we set $t= c-r_g$. From above assertion we get (ii) and (iii).
After all, we set $\gamma= \max \{1, d_0(a,b,c) \}$.
\end{proof}

\section{Almost-K\"{a}hler island metrics on $\Bbb R^{4}$}

In this section we describe almost \k  island metrics on $\Bbb
R^{4}$. We consider a metric on $\Bbb R^4$ of the form
\begin{equation} \label{g0fh}
 \tilde{g} = f^2 dr^2 + {r^2 \over f^2} d\theta ^2 + h^2
d\rho ^2 + {\rho ^2 \over h^2} d\sigma ^2,
\end{equation}
 where $(r,
\theta ) , (\rho ,\sigma)$ are the polar coordinates for each
summand of $\Bbb R^4 := \Bbb R^2 \times \Bbb R^2$ respectively and
$f , h $ are smooth positive functions on $\Bbb R^{4}$, which are
functions of $r$ and $\rho$ only. Let $e_1 = {1 \over f}{\partial
\over \partial r}$, $e_2 = { f \over r}{\partial \over
\partial \theta}$, $e _3 = {1 \over h }{\partial \over \partial
\rho}$,
    $e_4 = {h \over \rho} {\partial \over \partial \sigma}$.
A smooth almost complex structure $J_{\tilde{g}}$  is defined by
$J_{\tilde{g}}(e_1) =e_2 , J_{\tilde{g}}(e_2)=-e_1 ,
J_{\tilde{g}}(e_3)=e_4 , J_{\tilde{g}}(e_4)=-e_3$. With the standard
symplectic structure $\omega$ on $\Bbb R^4$ the triple $ (\tilde{g},
\omega, J_{\tilde{g}}) $ is an almost-K\"{a}hler structure on $\Bbb
R^4$. Let $\omega_i$ be the dual co-frame field of $e_i$. Compute
 the connection 1-forms ${\omega}_{ij}$ w.r.t.
$\omega_i$: $d {\omega}_{i } =
 \sum_{j=1}^{4}{\omega}_{ij} \wedge {\omega}_{j}$, with $ {\omega}_{ij}=
 -{\omega}_{ji}$;
one may compute $2a_{ijk} = \langle d {\omega}_{k}, {\omega}_{i}
\wedge {\omega}_{j} \rangle_g -\langle d {\omega}_{i}, {\omega}_{j}
\wedge {\omega}_{k} \rangle_g -\langle d {\omega}_{j}, {\omega}_{k}
\wedge {\omega}_{i} \rangle_g$, where $ {\omega}_{ij} = \sum_{k=1}
a_{ijk} \omega_k  $.
 And use the formula $d\omega_{ij} - {\omega}_{ik}\wedge {\omega}_{kj} = \sum_{k<l} R_{ijkl} \omega_k \wedge \omega_l $. Some Riemannian
curvature components of $\tilde{g}$  are computed;
\begin{align*}
&R_{1212} = -{3f_r \over rf^3} + {3f_r ^2 \over f^4 }- {f_{rr} \over f^3} - { f_{\rho}^2 \over f^2 h^2 },\\
%&R_{1223} = {2f_{\rho} \over rhf^2} + {f_{r \rho} \over rf^2 }- {3f_{\rho}f_r  \over hf^3} - { f_{\rho}h_{r} \over f^2 h^2 },\\
&R_{1313} = {f_{\rho\rho} \over fh^2} - {f_{\rho} h_{\rho} \over fh^3 }+ {h_{rr} \over f^2 h} - { f_r h_r \over f^3 h },\\
&R_{1414} = -{h_{rr} \over hf^2} + {h_r f_r h + 2f h_r ^2 \over f^3 h^2 }+ {f_{\rho} \over \rho f h^2} - { f_{\rho} h_{\rho} \over f h^3 },\\
%&R_{1434} = -{h_{r\rho} \over fh^2} + {h_r f_{\rho} h + 3 h_r fh_{rho} \over f^2 h^3 }-{2h_r \over \rho f h^2},\\
&R_{2424} = -{h_r \over rf^2 h} + {h_r f_r  \over f^3 h }- {f_{\rho} \over \rho f h^2} + { f_{\rho} h_{\rho} \over f h^3 },\\
%&R_{3414} = -{2h_r \over \rho h^2 f} - {h_{\rho r } \over f h^2 }+ {3 h_{\rho}h_r \over h^3 f} + { f_{\rho} h_r \over f^2 h^2 }\\
&R_{3434} = -{3h_{\rho} \over \rho h^3} - {h_{\rho \rho} \over h^3 }+ {3h_{\rho}^2 \over h^4} - { h_{r}^2  \over f^2 h^2 },\\
%&R_{2312} = {f_{\rho r} \over hf^2} - {3h_r f_{\rho} h + f_{\rho} h_r f \over f^3 h^2 }+ {2f_{\rho} \over rh f^2},\\
&R_{2323} = -{f_{\rho \rho } \over fh^2} + { f_{\rho} (2f_{\rho} h +
fh_{\rho}) \over f^2 h^3 }+ {h_r \over rf^2 h} - { f_r h_r \over f^3
h },
\end{align*}
where $f_{r}={{\partial f} \over {\partial r}},
f_{rr}={{\partial^{2} f} \over {\partial r \partial r}}$, etc.. The
scalar curvature is then as follows;
$$
{s_{\tilde{g}}}= - \{ (f^{-2})_{rr} +{3 \over r}(f^{-2} )_{r} +
(h^{-2})_{\rho \rho} + {3 \over \rho} (h^{-2})_{\rho} \}
                -{ 2f_{\rho}^2 \over h^2 f^2} - { 2h_{r}^2 \over h^2 f^2}.
$$
Setting $F = f^{-2}$ and $H = h^{-2}$, we shall find  $F$ and $H$
which satisfy
\begin{equation} \label{g1fh} F_{rr} + {3 \over r}F_r
+  H_{\rho \rho} + {3 \over \rho }H_{\rho}  = 0.
\end{equation}
\par
\noindent Set $F_{rr} + {3 \over r}F_r = \alpha (r) \beta (\rho)$
where $\alpha$, $\beta$ are smooth functions
on $\Bbb R$ %of $r$ and $\rho$ respectively
which satisfy at least that
\begin{align*}
 &\alpha(r) = 0  \quad for \quad r\leq 0 ,\ r \geq 1 ,\\
 &\beta(\rho)  = 0  \quad for \quad \rho\leq 0 ,\ \rho \geq 1.
\end{align*}
\par
\noindent These two functions $\alpha$ and $\beta$  will be
specified more below.

\medskip
Since $(r^3 F_r )_r = r^3 F_{rr} + 3r^2 F_r = r^3 \alpha (r) \beta
(\rho)$, we do integration with proper boundary conditions on $F$ to
get
%$F= \beta (\rho) \int ({1 \over r^3} \int r^3 \alpha (r) dr )dr$.
%We define the function $F$ to be
$$F(r, \rho ) = \beta (\rho) \int_{0}^r ({1 \over y^3} \int_{0}^y  x^3 \alpha (x) \,dx )\,dy +1. $$
\par
\noindent We now specify the function $\alpha$ as follows; first
consider a smooth function $p(y)$ on $\Bbb R$ such that
$$
\begin{cases}
&a)\  p(y) = 0 \quad \mbox{ for } \quad y\leq 0 , y\geq 1 ,\\
&b)\  |\frac{p'(y)}{y^3}| \ll 1 , \mbox{ for } y>0\\
&c)\  \int_{0}^1 {p(y) \over y^3} \,dy =0 , \\
&d)\  0<\int_{0}^r {p(y) \over y^3} \,dy <1 \ \ \mbox{ for  any } \
r \ \mbox{ with } \ 0<r<1
\end{cases}
%\leqno (?)
$$
and then define $ \alpha (y) = { p^{'} (y) \over y^3} $. Here we may
choose one such function $p$ so that the derivative $\alpha^{'}$ of
$\alpha$ is zero at exactly three points in the open interval
$(0,1)$ as in Fig.1.

\begin{xy}
 0="a",
(100,0)="b", "a";"b"**\dir{-}, (0,-20);(0,20)**\dir{-},
"a";"b"**\crv{(0,0)&(7,0)&(15,10)&(40,-10)&(60,-13)&(85,17)&(95,0)},
(-5,-15)*\txt{$-0.2$}, (-1,-15);(1,-15)**\dir{-},
(-5,15)*\txt{$0.2$}, (-1,15);(1,15)**\dir{-}, (15,-3)*\txt{$r_{1}$},
(15,-1);(15,1)**\dir{-}, (50,-3)*\txt{$r_{2}$},
(50,-1);(50,1)**\dir{-}, (85,-3)*\txt{$r_{3}$},
(85,-1);(85,1)**\dir{-}, (100,-3)*\txt{$1$},
(100,-1);(100,1)**\dir{-}, (-10,0);(0,0)**\dir{-},
(100,0);(110,0)**\dir{-}, (-1,-2)*\txt{$0$}
\end{xy}
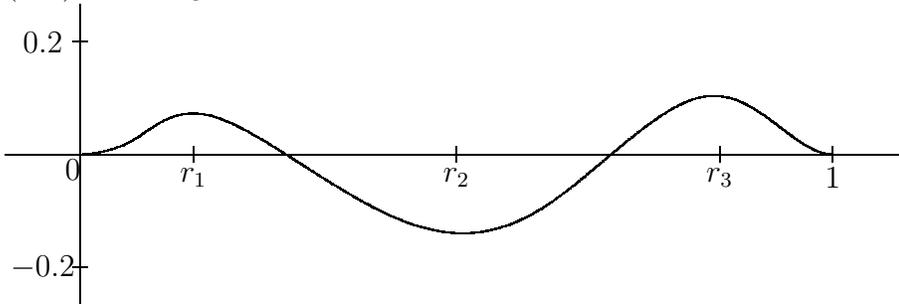
\centerline{Fig.1. The graph of $\alpha$.}
\medskip

\noindent Similarly we set
$$H(r,\rho) = -\alpha (r) \int_{0}^{\rho} ({1 \over y^3} \int_{0}^y  x^3 \beta (x) \,dx )\,dy +1 ,$$
%which satisfy $H_{\rho \rho} + {3H_{\rho} \over \rho }  = - \alpha (r) \beta (\rho)$.
where $\beta (y) = {k'(y) \over y^3}$ and the function $k(y)$ is a
smooth function on $\Bbb R$ satisfying $a)-d)$. %(\ref{g2fh}).
 We may choose
$k$ similarly to $p$ (or equally if one prefers) so that the
derivative $\beta^{'}$ of $\beta$ is zero at three points in
$(0,1)$. Hence the graph of $\beta$ is similar to that of $\alpha$.

\par
\noindent Then the functions $F(r,\rho)$ and $H(r,\rho)$ satisfy the
equation (\ref{g1fh}) and
\begin{align*}
&F,\ H \equiv 1 \ \ \ \mbox{for} \ r \leq 0 , \  r \geq 1 \ \ \mbox{or} \ \  \rho \leq 0 , \   \rho \geq 1  , \\
&F,\ H >0 \ \ \ \mbox{for} \ \ 0 < r < 1 \ \ \mbox{and} \ \ 0 < \rho
< 1 .
\end{align*}
In particular, the metric $\tilde{g}$ is smooth on $\Bbb{R}^4$. For
such functions $F$ and $H$, we have
\begin{align*}
s_{\tilde{g}} = & -  {H \over 2F^2 }  \beta ' (\rho ) ^2 \{\int_{0}^r ({1 \over y^3} \int_{0}^y  x^3 \alpha (x) \,dx )\,dy\}^2 \\
              &  -{F \over 2H^2 } \alpha ' (r) ^2 \{\int_{0}^{\rho} ({1 \over y^3} \int_{0}^y  x^3 \beta (x) \,dx )\,dy\}^2.
\end{align*}
From the choices of $\alpha$ and $\beta$ , the derivatives
$\alpha^{'}(r)$ and $\beta^{'}(\rho)$ are zero at three points $r_1
, r_2 , r_3 $ and $\rho_1 , \rho_2 , \rho_3$, respectively. It is
now simple to check

\begin{equation} \label{s12}
 s_{\tilde{g}}(r,\rho)
             =\begin{cases}
           0 & \rm{if} \ r=0,\ r \geq 1, \  \rho=0 \ \ or \ \rho \geq 1 ,\\
           0 & \rm{at} \ (r_{i}, \rho_{j}), \ i,j=1,2,3, \\
           \rm{negative} & \rm{if} \ \ r, \rho \in (0,1) \ \ with \ r \neq r_{i} \ or \ \rho \neq \rho_{j}.
           \end{cases}
\end{equation}

 \noindent Now $\tilde{g}$ is an island metric and as such
one may directly apply the main deformation of section 7 later, but
to make our presentation clear we modify  $\tilde{g}$ a little here.
First, by choice of  ingredient functions $\alpha$ and $\beta$, we
can assume that  $\tilde{g}$ satisfies $\|\tilde{g}- g_{0}\|_{{C^{2
\kappa + 4}_{g_{0}}} (\Bbb{R}^4)} \leq \frac{\lambda_3}{8}$, where
$g_0$ is the Euclidean metric $dr^2 + {r^2 } d\theta ^2 + d\rho ^2 +
{\rho ^2} d\sigma ^2$.

  We shall apply
the deformation of section 4 to $\tilde{g}$. Choose a point $p$ in
$\Bbb{R}^4$ such that $r(p)=\rho(p)= \frac{r_1}{50} $. Then
$s(\tilde{g})$
 is negative on $B^{\tilde{g}}_{\frac{r_1}{100}}(p)$ by (\ref{dist}).
In applying Lemma \ref{5} and Lemma \ref{6}, we let $c=1$, $a=0.01$,
$b= 1-\frac{r_1}{200}$ and $\epsilon= \frac{r_1}{600}$, i.e. we
consider the metric $\tilde{g}_{d,s}^{1-\frac{r_1}{200},
\frac{r_1}{600}, 1}$ defined in terms of the $\tilde{g}$-orthonormal
co-frame $\omega_i$ as defined in \rm{(\ref{oms})} with $ \omega_1 = dr_{\tilde{g}}$ and
 $r_{\tilde{g}}(x) = {\rm dist}_{\tilde{g}} (p, x)$.
 \noindent By
Lemma \ref{6}, for $d =d_0=  \gamma(0.01, 1-\frac{r_1}{200},1)$ and
$s \in [0,1]$ we have

\rm{(i)} $ \tilde{g}_{d,s}^{1-\frac{r_1}{200}, \frac{r_1}{600}, 1}
\equiv  \tilde{g}$ on $B^{\tilde{g}}_{ \frac{r_1}{300}}(p)$,
%$ \Bbb{R} \setminus (0, b+ \epsilon) \times \Bbb{R}^n $.

\rm{(ii)} $s(\tilde{g}_{d,s}^{1-\frac{r_1}{200}, \frac{r_1}{600},
1})-s(\tilde{g}) \leq 0$ on $\Bbb{R}^4 \setminus
B^{\tilde{g}}_{\frac{r_1}{200}}(p) $,
%$(0,b] \times \Bbb{R}^n$,

\rm{(iii)} $s(\tilde{g}_{d,s}^{1-\frac{r_1}{200}, \frac{r_1}{600},
1})-s(\tilde{g}) \leq -s \cdot e^{- \frac{d}{a}}$ on $
B^{\tilde{g}}_{0.99}(p) \setminus B^{\tilde{g}}_{\frac{r_1}{200}}(p)
$.

 From Lemma \ref{5} (i)
%for each $\varepsilon \in (0,1)$
we have  ${s}_0 $ with $1>{s}_0 >0 $ such that for $d
=d_0=\gamma(0.01, 1-\frac{r_1}{200},1), s^{-1} = {s}_0^{-1}$
%and each subset $B_i \subset A$:
we have
\begin{equation} \label{albe}
s(\tilde{g}_{d_0,s_0}^{1-\frac{r_1}{200}, \frac{r_1}{600}, 1}) < 0,
\mbox{ on } B^{\tilde{g}}_{\frac{r_1}{150}}(p)  \  \  \  \mbox{ and,
}
\end{equation}
\begin{equation} \label{albe}
\ \ \|\tilde{g}_{d_0,s_0}^{1-\frac{r_1}{200}, \frac{r_1}{600},
1}-\tilde{g}\|_{{C^{2 \kappa + 2}_{g_{0}}} ( \Bbb{R}^4)} <
\frac{\lambda_3}{8}.
\end{equation}

We shall denote by $g^-$ the metric
$\tilde{g}_{d_0,s_0}^{1-\frac{r_1}{200}, \frac{r_1}{600}, 1}$. As
$0<r_1 < 1$, from above we deduce the following proposition;

\begin{prop} \label{p2}
%Let $B^1_r\subset and $B \subset \Bbb R^{2n} $, $n\geq 3$, be the $n$-product of $B^{1}_r$.
With the standard symplectic structure on $\Bbb R^{4}$ there exists
a compatible almost-K\"{a}hler metric  $g^-$ of non-positive scalar
curvature on $\Bbb R^{4}$ which has the scalar curvature negative on
the ball $\{x| \ |x|< 0.9  \}$
%the polydisc $B$, the product of two open balls of radius  1  in $\Bbb R^2$,
and is Euclidean on  $\{x| \ |x|> 1.5  \}$. Furthermore, it
satisfies
 \begin{equation} \label {p2a}
 \| g^- -
g_{0}\|_{{C^{2 \kappa + 2}_{g_{0}}}(\Bbb{R}^4)} \leq
\frac{\lambda_3}{4}.
\end{equation}
\end{prop}

We remark that the metric in (\ref{g0fh}) was first exhibited in
\cite{KK}. Its scalar curvature, however, was not well presented,
which should be as in above (\ref{s12}). So for completeness sake we
rewrite it and apply the scalar-curvature diffusion, rather than
supplying extra to the argument of that paper.

\section{Almost \k Surgery}

Given a compact almost \k manifold $(M, g, \omega, J)$ of real
dimension $4$, for each point $p \in M$ there exists an open
neighborhood $V_p$ in which there exists a coordinate system
%$x_i$ $i=1, 2, \cdots n$, i.e. $\omega = \sum_{i=1}^{n} dx_i \wedge dx_{i+1}$.
where $\omega$ has the standard form. Henceforth we shall call such a
one  a {\it Darboux} coordinate system.
%Choose an open $g$-geodesic ball in  $ U_p$.
 We choose an open coordinate ball of finite radius $U_p$ such that $U_p \subset \overline{U_p} \subset V_p$,
 where $\overline{U_p}$ is the compact closure of $ U_p$ in $M$.
These $U_p$'s  form a covering of $M$ and there is a finite
subcovering $\{U_{p_{\alpha}}\}_{\alpha}$, which we denote by ${\frak A}$.
%, together with a Darboux coordinate system $y_i$ each.
We shall work with this (fixed) atlas.
 Denote  the Lebesgue number of ${\frak A}$ with respect to $g$  by $\sigma := \sigma_{g, \frak A}$.

\smallskip
 Then, near each point $p$ in $M$
 %of the above almost \k structure $(\omega,  g, J)$,
 there is a Darboux coordinate system $y_i$ in $\frak A$
 so that $\omega = \sum_{i=1}^{2} dy_{2i-1} \wedge
dy_{2i}$ on $B^{g}_{\sigma/2} (p)$. Then $ g$ can be written as $ g
= \sum_{i,j=1}^{4} g_{ij} dy_i \otimes dy_{j}.$ Consider the
Euclidean metric $ g_p= \sum_{i,j=1}^{4} g_{ij}(y(p)) dy_i \otimes
dy_{j} $ on $B^{g}_{\sigma/2} (p)$, i.e. we extend $g|_p$, the
restriction of $ g$ at $p$, onto $B^{g}_{\sigma/2} (p)$ with
constant coefficients. Because $y_i$ are Darboux coordinates, $g_p$
is $\omega$-almost \k with corresponding almost complex structure
$J_p$, the extension of $J|_p$ onto $B^{g}_{\sigma/2} (p)$ with
constant coefficients.
 Note that the definition of $g_p$ depends on the choice
of a Darboux coordinates system near $p$ in ${\frak A}$. However,
there are only a finite number of coordinate systems in ${\frak A}$,
defined with the above property. Due to this aspect the
following approximation Lemma holds, which provides the technical
basis for subsequent Lemmas.
%Here the rescaled metric $m^2 g_p$ is defined at least on $B^{m^2g}_{200} (p)$.

\begin{lemma} \label{fae}
The following two statements hold.
\par
\rm{(i)} For any $\varepsilon < 10^{-4}$, there is $m^{(0)} \geq
\frac{10^5}{\sigma}
>0$  with $m^{(0)}= m^{(0)}(\varepsilon)$ such that  for any $m >m^{(0)}$,
 $\|m^2 g - m^2 g_p\|_{{C^{2 \kappa + 2}_{m^2 g_p}} } <
\varepsilon$ on $B^{m^2g}_{10^4 + 1} (p)$ as well as on
$B^{m^2g_p}_{10^4} (p)$, where $p$ is any point of $M$.

\rm{(ii)} For any $\varepsilon < 10^{-4}$, there is $m^{(1)} \geq
\frac{10^5}{\sigma}
>0$  with $m^{(1)}= m^{(1)}(\varepsilon)$ such that for any $m >m^{(1)} $,
$\|m^2 g_q - m^2 g_p\|_{{C^{2 \kappa + 2}_{m^2 g_p}} } <
\varepsilon$ on $B^{m^2g}_{10^4 + 1} (p)$ as well as on
$B^{m^2g_p}_{10^4} (p)$, where $p$ is any point of $M$ and $ q \in
B^{m^2g}_{2 \cdot 10^4} (p) $.
\end{lemma}

\begin{proof}
Note that in the statement (ii) $m^2 g_q$ is defined on
$B^{m^2g}_{10^4 + 1} (p)$ because with $m  \geq
\frac{10^5}{\sigma}$, $m^2 g_q$ is defined on $ B^{m^2g}_{5
\cdot 10^4} (q)$.

 Let $e_i$,
$i=1, 2, 3,4$ be a parallel orthonormal frame field of $g_p$. Then
$e_i= \sum_{j=1}^{4} a_{ij} \frac{\partial}{\partial y_j}$, where
$a_{ij}$ are constants and $y_j$ is any Darboux coordinates used to define
$g_p$. By definition, in order to estimate $\|m^2 g
- m^2 g_p\|_{{C^0_{m^2 g_p}} }$ on $B^{m^2g}_{10^4 + 1} (p)$, we
 consider  $\mid m^2 g (\frac{1}{m} e_i, \frac{1}{m}e_j) -
m^2 g_p (\frac{1}{m}e_i, \frac{1}{m}e_j) \mid = \mid  \sum_{s,t}
a_{is} a_{jt} (g_{st}- {g_p}_{st}) \mid$
 on $B^{g}_{(10^4 +1)/m} (p)$.

 Note that $\max_{i,j} |{a_{ij}}|$ is determined by the continuous functions
$g_{st}= g(\frac{\partial}{\partial y_s }, \frac{\partial}{\partial
y_t} )$ defined on a coordinate ball $U_{p_{\alpha}}$ of finite
radius (there are only a finite number of such balls  in ${\frak
A}$) and choice of an orthonormal frame.
 So there is a bound $M$, i.e.
$|{a_{ij}}|
 \leq M$.
 %So $\mid  a_{is} a_{jt} \mid \leq L$, where $L$ is independent of $m$.
 For the same reason,
    for any $\varepsilon > 0$ there is $m^{(0)}$ such that
   $\mid \sum_{s,t} a_{is} a_{jt} (g_{st}- {g_p}_{st}) \mid \leq  M^2 | \sum_{s,t} (g_{st}- {g_p}_{st})|\leq  \varepsilon$
 on $B^{g}_{(10^4 +1)/m} (p)$   for $m >m^{(0)}$.
The $C^k$-case ($k \geq 1$) can be done similarly.
 This proves (i) on $B^{m^2g}_{10^4 + 1} (p)$.

Next, we claim that  $B^{m^2g_p}_{10^4} (p) \subset B^{m^2g}_{10^4
+1} (p)$; consider the $m^2g_p$-geodesic $l(t)$ with $l(0)=p$ and
$|l^{\prime}(t)|_{m^2g_p} \equiv 1$. Let $t_1$ be the first value of
$t$ such that $l(t)$ lies in the boundary of $B^{m^2g}_{10^4 +1}
(p)$. Then
$$ 10^4 +1 = d_{m^2 g}(p, l(t_1)) \leq \int_0^{t_1} | l^{\prime}(s) |_{m^2 g}
ds \leq \sqrt{1+ \varepsilon} \int_0^{t_1} |l^{\prime}(s) |_{m^2
g_p} ds \leq \sqrt{1+ \varepsilon} t_1. $$ As $10^4 < \frac{10^4 +
1}{\sqrt{1+ \varepsilon}}$, this proves $B^{m^2g_p}_{10^4} (p)
\subset B^{m^2g}_{10^4 + 1} (p)$, and finishes the proof of (i).

 The proof of (ii) is similar, but note that $g_p$ and $g_q$ may use different Darboux
 coordinates, say $y_i$ and $z_i$, respectively. We may write $g= \sum_{i,j} g_{ij} dy_i \otimes dy_{j}
 = \sum_{k,l} h_{kl} dz_k \otimes dz_{l} $,
  where $g_{ij} =\sum_{k,l} h_{kl} \frac{\partial z_k}{\partial y_i} \frac{\partial z_l}{\partial
  y_j}$.
For the $C^0$-norm $\|m^2 g_q - m^2 g_p\|_{{C^0_{m^2 g_p}} }$, it
suffices to
 show that  for $m$  large

\noindent
 $|m^2 g_q (\frac{e_i}{m}, \frac{e_j}{m})- m^2 g_p(\frac{e_i}{m}, \frac{e_j}{m})| \leq \sum_{s,t} \mid a_{is} a_{jt}
\mid \cdot \mid ({g_p}_{st} - g_{st}) +  \sum_{k,l} (h_{kl}-
h_{kl}(q)) \frac{\partial z_k}{\partial y_s} \frac{\partial
z_l}{\partial y_t} \mid \ \ $ is small on $B^{g}_{(10^4 + 1)/m} (p)$
for $ q \in B^{g}_{2 \cdot 10^4/m} (p) $. This can be done by a
similar argument to (i). The $C^k$-case ($k \geq 1$) can be done
similarly.
\end{proof}

For the given almost \k metric $g$ on $M$, by Lemma \ref{fae} and
Lemma \ref{lamb} if $m \geq m^{(0)}(\lambda_1)$, the injectivity
radius $\rm{inj}_{m^2g} \geq 200$ and the sectional curvature
$|K_{m^2g}| \ll 1$.

We recall the Besicovitch covering Lemma \ref{7}.
Assuming $m \geq m_0  \cdot m^{(0)}(\lambda_1)$, we shall explain a
surgery replacing neighborhoods of $p \in A=A(m)$ in the almost \k
manifold $(M,  m^2 g, m^2 \omega, J)$ with
 copies of the metrics $g^-$ of
Proposition \ref{p2}. Let $B_i :=\{ z \in A | \overline{B^{m^2
g}_{5} (z)} \in {\cal B}_i \}$. Then $A:= \cup_{i=1}^{\kappa} B_i$.

At each point $p \in A$, there is a Darboux coordinate system $y_i$
in ${\frak A}$ such that on $B^{m^2 g}_{10^4} (p)$, $m^2 \omega =m^2
\sum_{i} dy_{2i-1} \wedge dy_{2i}$, $m^2 g = \sum_{i} m^2 g_{ij}
dy_i dy^j$. And $m^2 g_p= \sum_{i} m^2 g_{ij}(p) dy_i dy_j$ is a
Euclidean metric compatible with $m^2 \omega$ and the corresponding
almost complex structure $J_p$. Now choose a $ g_p$-parallel
orthonormal frame field, $e_1$, $e_2= J_p e_1$, $e_3$,  $e_4=J_p
e_3$. There exists a coordinate
 system $z_i$ such that $e_i =\frac{\partial}{\partial z_i}$. In this coordinate $z_i$,
 $m^2 \omega =m^2  \sum_{i} dz_{2i-1} \wedge dz_{2i}$ and
 $m^2 g_p= \sum_{i} m^2 dz_i \otimes dz_i$.
Recall the metric $g^-$ of Proposition \ref{p2}
 on ${\Bbb R}^4= \{ x=(x_1,
x_2, x_3, x_4) \}$. Consider a diffeomorphism  $\phi = \phi_p :  \{
z | \mbox{ } |z| < \frac{10^4}{m} \} \rightarrow   \{ x | \mbox{
}|x|< 10^4 \}$ defined by $z \mapsto mz$. Then from the definition
of $g^-$ we note that  $\phi^{*} g^- = m^2 g_p$ outside $ \{ z | |z|
< \frac{1.5}{m} \}$.

\medskip From (\ref{geh}),  we may express $m^2 g =  (\phi^{*}g^- ) \cdot e^{h}$ on
$B^{m^2g}_{10^4} (p) $ for a unique smooth symmetric
$J_{\phi^{*}g^-}$-anti-invariant tensor $h$  because $m^2 g$ and
$\phi^{*}g^-$ are both $m^2 \omega$-compatible. Let
 $\eta (r)$  be a
smooth cutoff function in $C^{\infty}(\Bbb{R}^{>0}, [0,1])$ s.t.
$\eta \equiv 0$ for $0 < r<  1.7$ and $\eta \equiv 1$ on the set $ r
\geq 1.8$.

For $ m \geq m_0  \cdot m^{(0)}(\lambda_1) $ we define a metric on
$M$;
$$
g_A :=
\begin{cases}
& (\phi^{*}g^- )  \rm{e}^{{\eta}(r_{m^2 g}) \cdot h},  \hspace{.2in} \rm{on }
\cup_{p \in A} B_{2}^{m^2g} (p), \\
& m^2  g,   \hspace{0.9in} \rm{elsewhere}.
\end{cases}
$$

So we have embedded copies of the metric $g^-$ into $(M,
m^2 g, m^2 \omega, J)$ near each $p \in A$.

Note that from Lemma \ref{fae} (i) and formula (\ref{dist}) we have
$B_{1.5}^{m^2g_p} (p) \subset B_{1.7}^{m^2g} (p).$

\section{Main Deformation}
 Now we apply the scalar curvature diffusion of section 4, i.e. we
 deform $g_A$ so that we decrease the scalar curvature
away from neighborhoods of $p \in A$, and absorb the counter effect
(scalar curvature increase) in  neighborhoods of $p$  where the
scalar curvature of $g_A$ is negative.

Using the Besicovitch covering ${\cal B}_1, \cdots , {\cal
B}_{\kappa}$ we shall define $\omega$-almost \k metrics $g(i)$ on
$M$ inductively on $i=0, 1, 2, \cdots, \kappa$.

Observe that from (\ref{p2a}) we have that for $m \geq m_0 \cdot
m^{(0)}(\lambda_1)$,

\begin{equation} \label{gm15}
 \|\phi_q^* g^- - m^2
g_{q}\|_{{C^{2 \kappa +2}_{m^2g_{q}}}} \leq \frac{\lambda_3}{4}
\hspace{0.2in} \mbox{ on } B_{10^4}^{m^2g_q}(q), \mbox{ for any } q
\in A.
\end{equation}

\noindent We set $g(0) = g_A$.
%, where the metric $g^-$ used in defining $g_A$ is assumed to satisfy (\ref{gminus}).
 %\underline{after} such choices of $f$ and $h$.
We use often the following Lemma of base-metric change, whose proof is elementary.

\begin{lemma} \label{pq}
Let $g_1$, $g_2$ be two Euclidean metrics on a domain $\Omega$ in
$\Bbb{R}^4$ such that  $ B=B_{10^4}^{g_1} (p) \subset \overline{B}
\subset \Omega$. For any $\delta >0$ there exists
$\varepsilon^{(0)}=\varepsilon^{(0)}( \delta)$ such that if $\|g_1 -
g_{2}\|_{C^l_{g_1}(B)} < \varepsilon^{(0)}$, then $\|
\eta\|_{C^l_{g_1}(\cal{O})}\leq (1+ \delta)\|  \eta
\|_{C^l_{g_2}(\cal{O})} $ for any symmetric \rm{(2,0)}-tensor field
$\eta$ and any integer $l$ with  $2 \leq l \leq 2\kappa +2$,
% with $\|  \eta \|_{C^4_{g_2}(\cal{O})} \leq \mu$
where $ \cal{O}$  is any open subset of $B$.
\end{lemma}

Now we shall make estimates on $g(0)$.
Recall that  on the region $B_{1.8}^{m^2g}(q) \setminus B_{1.7}^{m^2g}(q)$ the metric
$m^2 g $ and the Euclidean metric $ \phi_q^* g^-= m^2 g_q$ are connected by a cutoff function.
As $m^2g \rightarrow m^2 g_q$ when $m \rightarrow \infty$, the next Lemma should be natural.
\begin{lemma} \label{111}
 There exists  $\tilde{m}_0$ with $ \tilde{m}_0 > m_0 \cdot m^{(0)}(\lambda_1)$
 such that for $m > \tilde{m}_0$,
$ \|g(0) - m^2 g_{q} \|_{{C^{2 \kappa +2}_{m^2g_q}} } \leq
\frac{\lambda_3}{4} \hspace{0.1in}$ on $B_{1.8}^{m^2g}(q) \setminus
B_{1.7}^{m^2g}(q)$
 for each $ q \in A$.
\end{lemma}

\begin{proof}
On $U_q=B_{1.8}^{m^2g}(q) \setminus B_{1.7}^{m^2g}(q)$, we define
$h_m$ in $m^2g=(m^2 g_q)  e^{h_ m}$. Then setting $m=1$, $h_1$ is
defined by $g= g_q e^{h_1}$. One easily sees from definition
(\ref{a001}) that $(m^2 g_q)^{-1}{h_m}= (g_q)^{-1} h_1 $ and
$\frac{h_{m}}{m^2} =h_1$. Note that in general  $(m^2 g)  e^{h} \neq
m^2 (g e^{h}) $.

If $\{e_i\}$ is a $
 g_q$-parallel orthonormal frame field, then $\{ \frac{e_i}{m} \}$
 is an $m^2 g_q$-parallel orthonormal frame field. On $U_q=B_{\frac{1.8}{m}}^{g}(q)
\setminus B_{\frac{1.7}{m}}^{g}(q)$, we have

\noindent $m^2 g(\frac{e_i}{m}, \frac{e_j}{m})= ( m^2 g_q) e^{h_m}
(\frac{e_i}{m}, \frac{e_j}{m}) =
% m^2g_q(\frac{e_i}{m}, e^{(m^2 g_q)^{-1}{h_m}} (\frac{e_j}{m})) =
 g_q({e_i}, e^{( g_q)^{-1}{h_1}}
({e_j}))$
 approach $\delta_{ij}$ as $m \rightarrow \infty$ uniformly, i.e. independently of
 point $q$ and the coordinates chosen to define $g_q$.
Or, simply it is  Lemma \ref{fae} (i). Clearly $g_q \rightarrow g|_q$ and $( g_q)^{-1}{h_1}\rightarrow 0$ uniformly.

Then, as $g(0) (\frac{e_i}{m}, \frac{e_j}{m})= (m^2 g_q)
\rm{e}^{{\eta}(r_{m^2 g}) \cdot h_m} (\frac{e_i}{m}, \frac{e_j}{m})
=  g_q({e_i}, e^{\eta(r_{m^2g})( g_q)^{-1}{h_1}} ({e_j}))$, so
$g(0)$
%$m^2 g_q \rm{e}^{{\eta}(r_{m^2 g}) \cdot h_m}$
approaches $m^2 g_q e^0
= m^2 g_q$ in the norm $\| \cdot \|_{C^0_{m^2g_q}}$.
The proof of convergence in $\| \cdot
\|_{C^{2 \kappa +2}_{m^2g_q}}$  is similar and straightforward.
\end{proof}

The next Lemma is essential in the construction of
the metrics $g(i)$'s;
\begin{lemma} \label{100}
There exists $m_1$ with $m_1 \geq \tilde{m}_0$ such
that for $m \geq m_1$,  $ \|g(0) - m^2 g_p \|_{{C^{2 \kappa
+2}_{m^2g_p}} } \leq \frac{\lambda_3}{2}$ on $B_{10^4}^{m^2g_p}(p)$
for each $p \in A$.
\end{lemma}

\begin{proof}
 Fix $p \in A$.
 On $V_1= B_{10^4}^{m^2g_p}(p) \setminus \cup_{q \in A} B_{1.8}^{m^2g}(q)$, $g(0) =m^2 g$.
By Lemma \ref{fae} (i), if $m \geq m^{(0)}(\frac{\lambda_3}{2})$
then
 $ \|g(0) - m^2 g_p \|_{C^{2 \kappa +2}_{m^2g_p}(V_1)}  \leq \frac{\lambda_3}{2}$.

 On $ B_{10^4}^{m^2g_p}(p) \cap  B_{1.7}^{m^2g}(q) $ with ${q \in A}$, $g(0)
= \phi_q^* g^-$, so we have the estimate (\ref{gm15}) if $ m \geq
m_0 \cdot m^{(0)}(\lambda_1)$.
 In the transition region $ B_{10^4}^{m^2g_p}(p)\cap
  [ B_{1.8}^{m^2g}(q) \setminus B_{1.7}^{m^2g}(q)]$, $q \in A$,
  Lemma \ref{111} can be applied if $m > \tilde{m}_0$.
 So, on $W_q= B_{10^4}^{m^2g_p}(p)
\cap B_{1.8}^{m^2g}(q) $  for $q \in A$ we estimated $g(0) - m^2
g_{q}$  in terms of $ m^2g_q$.

Now it remains to estimate $g(0) - m^2 g_{p}$ on $W_q$ with respect
to $m^2 g_{p}$.
  Applying Lemma \ref{fae} (ii) and Lemma \ref{pq} (with $g_1=m^2 g_p $ and $g_2= m^2 g_q$); if
$m \geq m^{(1)} (\varepsilon^{(0)}(\frac{1}{2}))$;
\begin{align}
\|g(0) - m^2 g_{p}\|_{C^{2 \kappa +2}_{m^2g_p}(W_q)}
 & \leq \|g(0) - m^2 g_{q} \|_{{C^{2 \kappa +2}_{m^2g_p}} (W_q)}
+  \| m^2g_q- m^2 g_p \|_{{C^{2 \kappa +2}_{m^2g_p}} (W_q)} \nonumber \\
&  \leq \frac{3}{2} \|g(0) - m^2 g_q \|_{{C^{2 \kappa +2}_{m^2g_q}}
(W_q)}
 +   \| m^2g_q- m^2 g_p \|_{{C^{2 \kappa +2}_{m^2g_p}} (W_q)}  \nonumber
\end{align}

\noindent By Lemma \ref{fae} (ii), $ \| m^2g_q- m^2 g_p \|_{{C^{2
\kappa +2}_{m^2g_p}} (W_q)} \leq \frac{\lambda_3}{10 \kappa} $
 if
$m \geq m^{(1)} ( \frac{ \lambda_3}{10\kappa})$.

\noindent From (\ref{gm15}) and Lemma \ref{111}, for $m \geq m_1 = \max \{
\tilde{m}_0, \mbox{ } m^{(1)} (\varepsilon^{(0)}(\frac{1}{2})),
\mbox{ } m^{(1)} ( \frac{\lambda_3}{10 \kappa}) \}$ we have
\begin{equation}
\|g(0) - m^2 g_{p}\|_{C^{2 \kappa +2}_{m^2g_p}(W_q)} \leq
\frac{3}{2} \cdot \frac{\lambda_3}{4} + \frac{\lambda_3}{10 \kappa}
<  \frac{\lambda_3}{2} .
\end{equation}

\noindent We proved  that  $ \|g(0) - m^2 g_p \|_{{C^{2 \kappa
+2}_{m^2g_p}} (B_{10^4}^{m^2g_p}(p))} \leq \frac{\lambda_3}{2}$ for
each $p \in A$ if $m \geq m_1.$
\end{proof}

Suppose that for $m \geq m_1$, $g(j)$ is defined when $j=0, \cdots
,i-1$ and satisfies the hypothesis
\begin{equation} \label{a12}
\|g(j) - m^2 g_p\|_{{C^{2\kappa +2 -2 j}_{m^2g_p}}
(B_{10^4}^{m^2g_p}(p))} \leq \{ \frac{1}{2} + \frac{j}{10 \kappa} \}
\lambda_3 \hspace{0.2in} \mbox{ for } p \in A.
\end{equation}
Lemma \ref{100} corresponds to the case $j=0$. We note that in
(\ref{a12}) $C^{2\kappa +2 -2j}$-norm is used because in defining $
g(1), \cdots, g(\kappa)$ successively we lose two derivatives at
each step in applying Proposition \ref{lemom} and Lemma \ref{5}, as
we shall see in the proof of the next lemma.

 Now, by Proposition \ref{lemom} on each $B_{100}^{g(i-1)}
(p)$ for $p \in B_{i}$, the $g(i-1)$-orthonormal co-frame field
$\omega_1, \omega_2, \omega_3, \omega_4$ of (\ref{oms}) is defined.
%and satisfy Proposition \ref{lemdr} and Lemma \ref{lemaijk}.
We set as in (\ref{gbe});
 \begin{align}
g(i)  = & \frac{ 1}{( 1+ h_{\epsilon}^{b} \cdot f_{d,s}  (9 - d_{
g(i-1)} (p,
x)))}\omega_1 \otimes \omega_1 \nonumber \\
 & + ( 1+ h_{\epsilon}^{b} \cdot f_{d,s} (9
- d_{ g(i-1)} (p, x)))
 \cdot \omega_2 \otimes \omega_2
 & + \omega_3 \otimes \omega_3 +\omega_4 \otimes \omega_4.
 \end{align}

\noindent with $ b= 8.5, a=1,  \epsilon=0.1, c=9$. With this done on
$\cup_{p \in B_{i}} B_{100}^{g(i-1)}(p)$, $g(i)$ is smooth and
$\omega$-almost \k on $M$ because $B_{9}^{g(i-1)} (p)$'s, $p \in
B_{i}$ are disjoint by Lemma \ref{7} (ii), Lemma \ref{fae} (i), the
formula (\ref{dist}) and (\ref{a12}).
%$$g(B, m, r,d,s) := \Pi_{p \in B} H^{9.5}_{0.1} \cdot F_{d,s}(9 - d_{m^2 \cdot g} (p, x))) \cdot g(0)$$

\medskip
\begin{lemma} \label{101}
There exists $\tilde{d}$, independent of $i$, such that for $d \geq
\tilde{d}$, $s \in (0, 1]$ and $m \geq m_1$, the metric $g(i)$
satisfies the hypothesis \rm{(\ref{a12})}, i.e. \rm{(\ref{a12})}
holds when $j=i$.
\end{lemma}

\begin{proof}
Consider $p \in A$ and let $z \in B_{10^4}^{m^2g_p}(p)$. Set $V_z=
B_{10^4}^{m^2g_p}(p) \cap  B_{1}^{m^2g}(z)$.

Note that $V_z = V_z \cap (\cup^{\kappa}_{j=1} \cup_{q \in B_j}
B^{m^2 g}_{10} (q)) = \cup^{\kappa}_{j=1} \{ V_z \cap (\cup_{q \in
B_j} B^{m^2 g}_{10} (q))\}.$  $g(i) = g(i-1)$ away from $\cup_{q \in
B_i} B^{g(i-1)}_{9} (q) \subset \cup_{q \in B_i} B^{m^2 g}_{10}
(q)$. So, we get
 $\|g(i) - g(i-1)\|_{C^{2\kappa +2 -2 i}_{m^2 g_p} (V_z)} =
 \|g(i) - g(i-1)\|_{C^{2\kappa +2 -2i}_{m^2 g_p}
 (V_z \cap (\cup_{q \in B_i} B^{m^2 g}_{10} (q)) )}$.
 %By Bishop-Gromov volume comparison argument, the number of $q$'s in $B_i$ such that
 %$V_z \cap  B^{m^2 g}_{10} (q) \neq \phi$ is finite.
  Let $Q$ be the set of
 $q$'s in $B_i$ such that
 $V_z \cap  B^{m^2 g}_{10} (q) \neq \phi$.
 For $q \in Q$, from Lemma \ref{pq} and Lemma
\ref{fae} (ii), if $m \geq m^{(1)} (\varepsilon^{(0)}(\frac{1}{2}))
$,
\begin{align}
 \|g(i) - g(i-1)\|_{C^{2\kappa +2 -2 i}_{m^2 g_p} (B^{m^2
g}_{10} (q) )} &= \|g(i) - g(i-1)\|_{C^{2\kappa +2 -2 i}_{m^2 g_p}
(B^{g(i-1)}_{9} (q) )} \nonumber \\
& \leq \frac{3}{2} \|g(i) - g(i-1)\|_{C^{2\kappa +2 -2 i}_{m^2 g_q}
(B^{ g(i-1)}_{9} (q) )}. \nonumber
\end{align}
 And from the  hypothesis (\ref{a12}) for $j=i-1$ and Lemma \ref{5} (ii),
 there exists $\tilde{d}= d( \frac{2\lambda_3}{30\kappa}, 9)$
%, $j=0, 1, \cdots, i$
so that if $d \geq \tilde{d}$ and $s \in (0, 1]$, $ \|g(i) -
g(i-1)\|_{C^{2\kappa +2 -2 i}_{m^2 g_q} (B^{ g(i-1)}_{9} (q) )} \leq
\frac{2\lambda_3}{30\kappa} $.

\noindent Combining  above, we have
\begin{align}
\|g(i) &- g(i-1)\|_{C^{2\kappa +2 -2 i}_{m^2 g_p} (V_z  \cap
(\cup_{q \in B_i} B^{m^2 g}_{10} (q)))} \\
 \leq  &  \sup_{q \in
Q} \|g(i) - g(i-1)\|_{C^{2\kappa +2 -2 i}_{m^2 g_p} (B^{m^2 g}_{10}
(q) )}  \leq \frac{3}{2} \cdot \frac{2\lambda_3}{30\kappa}=
\frac{\lambda_3}{10 \kappa}. \nonumber
\end{align}

\begin{align}
\mbox{Finally, } \| g(i) &- m^2 g_p  \|_{C^{2\kappa +2 -2 i}_{m^2 g_p} (V_z)}    \nonumber \\
 & \leq  \|g(i) - g(i-1)\|_{C^{2\kappa +2 -2 i}_{m^2 g_p} (V_z)}
 + \|g(i-1) - m^2 g_p\|_{C^{2\kappa +2 -2 i}_{m^2 g_p} (V_z)} \nonumber \\
& \leq   \frac{\lambda_3}{10\kappa} + \{ \frac{1}{2} + \frac{i-1}{10 \kappa} \} \lambda_3 = (
\frac{1}{2} + \frac{i}{10 \kappa} ) \lambda_3. \nonumber
\end{align}
So we proved (\ref{a12}) for $j=i$ and this proves the lemma.
\end{proof}

\noindent From Lemma \ref{101},  all the metrics $g(0)$, $g(1),
\cdots,  g(\kappa)$ can be defined. The parameters $d$ and $s$ are
being used in defining each metric $g(i)$, $i=1, \cdots, \kappa$.
But we shall use the same values of $d,s$ independently of $i$. Now
we can prove;

\begin{prop}
There are $m_2$ with $m_2 > m_1$ and $d_1, s_1 >0$  such that for
each $d \geq d_1$, $s^{-1} \geq s^{-1}_{1}$ there exists an $m \geq
m_2$ with $m= m(d, s)$ such that $s(g(\kappa)) < -c_1$ for a
constant $ c_1= c_1(d,s) >0$.
\end{prop}

\begin{proof}
 For $\delta >0$, there is an $m^{(2)} > m_1 $ with $m^{(2)} =
m^{(2)}( \delta) $ and
\begin{equation} \label{abcd}
 \sup_M  |\mbox{ }s(m^2 \cdot
g) | < \delta \text{        and       } - \mu -\delta < s(g(0)) <
\delta,
\end{equation}
for  $m \geq m^{(2)}$ and $\mu := - \rm{min}
\{ s(g^-) \} >0$.
%Recall that $g^-$ satisfies (\ref{gminus}).
Here, for the second inequality one use the argument of Lemma
\ref{111}, i.e. $ g(0)$ approaches $m^2 g_{q}$  on the transition
region $B_{1.8}^{m^2g}(q) \setminus B_{1.7}^{m^2g}(q)$ as $m$
becomes large.

\par
 We consider $m > m_1$.
Recall from Proposition \ref{p2} that $s(g^-) < 0$ on $\{x| \ |x|<
0.9  \}$. So from Lemma \ref{101}, Lemma \ref{5} (i) and Lemma
\ref{6} (ii)
%for each $\varepsilon \in (0,1)$
we have $\bar{d}_0$ with $\bar{d}_0 \geq \rm{max}\{ \tilde{d}, \gamma(1, \ 8.5,
\ 9) \}$ and
$\bar{s}_0 $ with $1 \geq \bar{s}_0
>0 $ such that if we use $d, s$ with $d \geq \bar{d}_0, s^{-1} \geq \bar{s}_0^{-1}$
in defining each metric $g(i)$,
%and each subset $B_i \subset A$:
we have
\begin{equation} \label{albe}
%& ({\rm i}) \mbox{ } |g(c) -g(0)|_{C^4_{m^2 g}} < \varepsilon,  \nonumber \\
 s(g(\kappa)) < -\beta \  \mbox{ on } B^{m^2 g}_{0.8}(p), \ p
\in A
\end{equation}

%& ({\rm i}) \sup _{p \in B_{i+1}}|g(i+1) -g(i)|_{C^4_{m^2 g_p}} < \varepsilon, \nonumber  \\
\noindent where  $\beta >0$ is a constant.
%We calculate the scalar curvatures $s(g(i+1))$ from $s(g(i))$.

By Lemma \ref{6} for $ d \geq \bar{d}_0$ and $s^{-1} \geq
\bar{s}_0^{-1}$, \begin{equation} \label{ssd}
 s(g(i+1)) -s(g(i)) \leq
\begin{cases}
&  -s \cdot e^{-d}  \text{  on  } \cup_{p \in B_{i+1}}
B^{g(i)}_{8}(p) \setminus B^{g(i)}_{0.5} (p),
\\
& 0   \hfill \text{ on  } M \setminus \cup_{p \in B_{i+1}}
B^{g(i)}_{8}(p).
\end{cases}
\end{equation}

Note that the balls $B^{g(i)}_{8}(p)$, $p \in B_{i+1}$ are disjoint
as $m
> m_1$ from Lemma \ref{7}(ii), Lemma \ref{fae} (i), Lemma \ref{101} and the formula
(\ref{dist}).
 Now the rest of proof is almost the same as that in
\cite[proposition 4.5]{Lo1}. We produce it for completeness sake.

By Lemma \ref{fae} (i), Lemma \ref{101} and the formula (\ref{dist})
again, we have $B^{m^2 g}_{6} (p)
 \setminus B^{m^2 g}_{0.55} (p)  \subset B^{g(i)}_{8}(p) \setminus B^{g(i)}_{0.5} (p)$ for $p \in B^{i+1}$.
%By definition  $B^{m^2 g}_6(p)$'s, $p \in A$ cover $M$, so

Therefore on $M \setminus  \cup_{p \in A} B^{m^2 g}_{0.6} (p)
\subset \cup_{p \in A} B^{m^2 g}_{6} (p)
 \setminus B^{m^2 g}_{0.55} (p)$,
by adding inequalities (\ref{ssd}), $i=0 ,   \cdots , \kappa-1,$ we
get $ s(g(\kappa)) - s(g(0)) \leq -{s \cdot e^{-d}}.$ From
(\ref{abcd}) we have
$$ %-s \cdot \kappa  \cdot c_1  - \mu -\delta <
s(g(\kappa ))  < -{s \cdot e^{-d}}+ \delta  \ \  \mbox{ on } M
\setminus  \cup_{p \in A} B^{m^2 g}_{0.6} (p),$$ for any $\delta
>0$, $m \geq   m^{(2)} (\delta)$  and $d \geq
\bar{d}_0, s^{-1} \geq \bar{s}_0^{-1}$. Note that $\bar{d}_0,
\bar{s}_0$ are independent of $\delta$. Choose $ \delta = \frac{s
\cdot e^{-d}}{2}.$ Then from (\ref{albe}) we obtain on $M$,
$% -\rm{max} \{  s \cdot \kappa \cdot c_1  + \mu  +  \frac{s \cdot e^{-d}}{2}, \alpha    \}<
 s(g(\kappa))     < -\rm{min} \{ \frac{s \cdot e^{-d}}{2}, \beta \}.$

Hence we will choose  $d_1 = \bar{d}_0$, $s_1 =\bar{s}_0 $, $m_2 =
 m^{(2)}( \delta_0) $ with $ \delta_0 =
\frac{\bar{s_0} \cdot e^{-\bar{d_0}}}{2}$ and find that for each  $d
\geq d_1$, $s^{-1} \geq \bar{s}_1^{-1}$ there is an $m \geq m_2$
with
\begin{equation*} s(g(\kappa)) < -c_1 \hspace{0.2in} \mbox{ for
} \hspace{0.2in} c_1 = \rm{min} \{ \frac{s \cdot e^{-d}}{2}, \beta
\} .\end{equation*}
\end{proof}
\noindent We now have an almost \k structure $(g(\kappa), m^2
\omega)$ of negative scalar curvature. The re-scaled metric $
\frac{g(\kappa)}{m^2}$ is $\omega$-compatible with negative scalar
curvature.  Theorem \ref{th1} is now proved.

\begin{remark}
\rm{Note that a \k version of the main theorem cannot hold. Indeed,
we have
 the well known formula $ \int_M s(g) \mbox{ } dvol_g=  4\pi c_1(J)
\cdot [\omega],$ where  $c_1(J)$ is the first Chern class of a
compact \k surface $(M, g, \omega, J)$ and $[\omega]$ is the
cohomology class of $\omega$.  So, on the complex projective plane
$\Bbb{CP}_2$ or a Hirzebruch surface $F$, a negative-scalar-curved
\k metric does not exist, whereas negative-scalar-curved almost \k
metrics exist on them by our main theorem.}
\end{remark}

\begin{remark}
\rm{ In this article we only worked on the 4-dimensional case. And
we wish to study the higher dimensional case in a forthcoming paper.
It is interesting to pursue any implication of this result to
symplectic topology.}
\end{remark}

\begin{remark}
\rm{From the nature of the scalar-curvature deformation argument of
section 4,
 one may speculate that if there is a scalar-island metric in some category of metrics, then one may prove the
existence of negative scalar curvature in the category. For instance
one may study Hermitian metrics or contact metrics etc..}
\end{remark}

\begin{remark}
\rm{ It is interesting to know if every closed symplectic manifold
of dimension four admits a compatible almost \k metric of negative
{\it Ricci} curvature or if the space of almost \k metrics with
negative scalar curvature is contractible, as in \cite{Lo3}.}
\end{remark}
\bigskip

\appendix
{\hspace{5cm} \bf{Appendix}}

\medskip
In this appendix we show the proof of Proposition \ref{lemom} and
supply details for the scalar curvature formula (\ref{sgt}).

\medskip
By Lemma \ref{lamb}, if   we assume $\|g- g_{0}
\|_{{C^2_{g_{0}}}(B^{g_{0}}_{10^4}(p))} \leq \lambda_1$, then $g$
has ${\rm inj}_g(p) \geq 200$. Let $B$ denote the set $\{ v \in T_p
\Bbb{R}^4 \mid \ |v|_{g|_p} \leq 100 \}$, where $g|_p$ is the
restriction of $g$ to $T_p \Bbb{R}^4$. Then the $g$-exponential map
at $p$, ${\rm exp}_p: B \rightarrow B^{g}_{100}(p) \subset \Bbb{R}^4
$ has a smooth inverse ${(\rm exp}_p)^{-1}$.
%As $ B^{g}_{100}(p)$ is a subset of $\Bbb{R}^4$,
If we identify $T_p \Bbb{R}^4 $ with $\Bbb{R}^4$ via the the map $i:
T_p \Bbb{R}^4 \rightarrow \Bbb{R}^4$, $ i(\sum_{i} v_i
\frac{\partial }{\partial x_i}|_p )= (v_1, \cdots , v_4)$,
%On $T_p {\Bbb R}^4$, we use the coordinate induced by $\frac{\partial}{\partial x_i}\mid_p$.
then we can view ${\rm exp}_p$ as a map from a domain $i(B)$ in
$\Bbb{R}^4$ to $\Bbb{R}^4$. For such a map we shall use the $C^k$
norm $\| \cdot \|_{C^{k}}$ in the next lemma.

\begin{lemma} \label{106}
There exists $\lambda^{(0)}$  with $\lambda^{(0)} < \lambda_1$  and
a positive constant $C^{(0)}$ such that
%For any $\varepsilon >0$, there exists $\delta \leq \lambda_1$, depending only on $\varepsilon$,
 if $\|g- g_{0} \|_{{C^l_{g_{0}}}(B^{g_{0}}_{10^4}(p))} < \lambda \leq \lambda^{(0)}$ for $ 4 \leq l \leq 2\kappa +4$, then $\| {\rm
exp}_p - ({\rm Id} + p) \ \|_{C^{l-1}(i(B))} < C^{(0)} \cdot
{\lambda}$ and $\| {(\rm exp}_p)^{-1} - ({\rm Id} -p) \
\|_{C^{l-1}(B^{g}_{100}(p))} < C^{(0)} \cdot {\lambda}$. Here, ${\rm
Id} + p $ is the map $({\rm Id} + p ) (x) = x + p$, and likewise for
${\rm Id} - p $.
\end{lemma}

\begin{proof}
%First we compare the derivatives of $\omega_1= dr_g$ and $\omega^0_1= dr_{g_0}$.
%We use the standard coordinates $(x_1, x_2, x_3, x_4)$ in which $g_{0} =\sum_{i=1}^{4}  dx_i \otimes dx_{i}$
%and $ g= \sum_{i,j=1}^{4} g_{ij} dx_i \otimes dx_{j}$.
To deal with the
geodesic
 equation of the metric $g= \sum_{i,j=1}^{4} g_{ij} dx_i \otimes
 dx_{j}$, i.e.
$ \ddot{x}_k = -\sum_{i,j=1}^4 \Gamma_{ij}^k( x(t)) \dot{x}_i
\dot{x}_j, $ where $ \Gamma_{ij}^k$ is the Christoffel symbol of
$g$,
 we consider the equivalent first order differential
system:
\begin{equation} \left\{ \begin{array}{ll}
\dot{x}_k & = y_k  \\
\dot{y}_k & = - \sum_{i,j} \Gamma_{ij}^k( x(t)) y_i y_j
          \end{array}
          \right. \label{geod} \end{equation}

\noindent
%without loss of generality we may only consider geodesics emanating from the origin; we set
with the initial condition $x(0) = p$ and $y(0)=v$. Denote its
solution curve by $(x(t, v), y(t, v))= ( x_1(t, v), \cdots, x_4(t,
v), y_1(t, v), \cdots, y_4(t, v))$.
%Consider the maps ${\rm exp}_p$ and $({\rm exp}_p)^{-1}$ in  $T_p {\Bbb R}^4  \cong {\Bbb R}^4
%\stackrel{\rm{exp_p}}{\longrightarrow}   {\Bbb R}^4\stackrel{\rm{exp_p}^{-1}}{\longrightarrow {\Bbb R}^4, $ and
%On $T_p {\Bbb R}^4$, we use the coordinate induced by $\frac{\partial}{\partial x_i}\mid_p$.

Suppose that $\|g- g_{0} \|_{{C^l_{g_{0}}}} \leq \varepsilon \leq
\lambda_1 < 10^{-4}$ with $ 4 \leq l \leq 2\kappa +4$.
%on $B^{g_{0}}_{11}(0) \subset {\Bbb R}^n$
Then, we have $|g_{ij}- \delta_{ij} | \leq \varepsilon$,
$|\frac{\partial^r g_{ij} }{\partial x_{m_1}  \cdots
 \partial x_{m_r}}| \leq \varepsilon$ for $r = 1,2,  \cdots, l$.  For $\Gamma_{ij}^k =
\frac{1}{2} \sum_{r=1}^4 g^{rk} ( \frac{\partial g_{jr}}{\partial
x_i} + \frac{\partial g_{ir}}{\partial x_j} -     \frac{\partial
g_{ij}}{\partial x_r} ) $, we have for $ s = 0, 1, \cdots, l-1$,
\begin{equation} \label{inv1}
|\frac{\partial^s \Gamma_{ij}^k }{\partial x_{m_1}  \cdots
 \partial x_{m_s}}| \leq  c_1 \varepsilon, \mbox{ with a positive
 constant } c_1.
 \end{equation}
\noindent In this appendix we shall use $c_i$, $i=1, 2, \cdots$ to
denote positive constants.
%Then $|\dot{y}| \leq | \frac{ \sum_i y_i \cdot \dot{y}_i}{y} | \leq 4 \varepsilon_1 y^2= 4 c_2 \varepsilon y^2$.
%the $g$-distance of the point $x(t, v)$ from $p$.
%As $\|g- g_{0} \|_{{C^4_{g_{0}}}} \leq \varepsilon< 10^{-4}$, it is enough to
%consider where $t |v| < 101$. So, setting

Since we study $ {\rm exp}_p$ on $B$, it is enough to consider $v$
with $ |v|:=|v|_{g_0} < 101$. And as $ {\rm exp}_p(v) =x(1, v)$, we
only need to consider $t$ with $ 0 \leq t \leq 1$ when dealing with
the solutions of (\ref{geod}). We write for convenience $ {\rm
exp}_p(v) = f(v)= (f_1(v), \cdots, f_4(v))$ and ${\rm exp}_p^{-1}(x)
= h(x) = (h_1(x), \cdots, h_4(x))$.

 \noindent Setting
$Y(t)=|y(t)|= \sqrt{\sum_{i=1}^4 (y_i(t) )^2}$, we get $ |\dot{Y}|
\leq c_2 \varepsilon \mbox{ } Y^2 $ from (\ref{geod}).
%Then we choose $\varepsilon$ so that $\varepsilon< \frac{1}{202 c_2} < \frac{1}{ 2 c_2 t |v|} $ and
Integrating it, we deduce from (\ref{geod}) $ \
%\frac{|v|}{1+ c_2 \varepsilon |v| t  }  \leq y(t) \leq \frac{|v|}{1-    c_2 \varepsilon |v| t }  \mbox{  } \mbox{  and  } \mbox{  }
|\dot{y}_k| \leq  c_2 \varepsilon Y^2 \leq \frac{c_2 \varepsilon
|v|^2}{(1-  c_2 \varepsilon t |v|)^2 } \leq c_3 \varepsilon$, for
small $\varepsilon$. So we get
\begin{align} \label{idon}
\hspace{.1in} & |y(t,v) - v|
%= | \int_0^t \dot{y}^k (s) ds|
\leq c_4 \varepsilon \mbox{,} \hspace{.4in} \
|x(t,v) -p -vt | \leq c_4 \varepsilon    \hspace{.1in} \mbox{} \hspace{.3in}\nonumber \\
& |f(v) - p - v| \leq c_4 \varepsilon  \mbox{,} \hspace{.4in}  | h
(x) - (x-p)| \leq c_5  \varepsilon .
\end{align}
%For $x \in B_{100}^g (p)$, we write $h(x) = v$ and $ | h (x) - (x-p)| =|v-(f(v)-p)|=|f(v)-p-v|< c_3
%\varepsilon |v|= c_3 \varepsilon |h(x)|$. For $\varepsilon $ small, we get $|h(x)| \leq c_4 |x-p|$.So,
where the last inequality is deduced from the others and $h \circ f
(v) = v$.

 Now, the solution $ x(t, v)$ and   $y(t, v)$ of
 (\ref{geod}) are well known to be $C^{\infty}$ functions of both variables $t$ and
$v$. Let $p_k$ and $v_k$ be the $k$-th component of $p$ and $v$,
respectively. From (\ref{geod}) we get for any $n$ with $1 \leq
n\leq l-1$,
\begin{equation} \left\{ \begin{array}{ll}
\frac{\partial^n x_k}{\partial v_{i_1}  \cdots \partial v_{i_n}} &
=\int_0^t \frac{\partial^n y_k}{\partial v_{i_1}  \cdots \partial
v_{i_n}} \ ds \\
 \frac{\partial^n (y_k - v_k)}{\partial v_{i_1}  \cdots \partial
v_{i_n}} & =  \int_0^t  \sum_{m,i,j} - \frac{\partial
\Gamma_{ij}^k}{\partial x^m} ( \frac{\partial^n x_m}{\partial
v_{i_1} \cdots \partial v_{i_n}}  ) y_i y_j - 2 \Gamma_{ij}^k (
 \frac{\partial^n y_i}{\partial v_{i_1}  \cdots \partial v_{i_n}} ) y_j + L_n ds
          \end{array}
          \right. \label{dyv} \end{equation}
\noindent where $L_n$ consists of terms involving partial
derivatives of $x$ and $y$ of order lower than $n$ only. Note that
$L_1=0$. We shall prove the following for $ |v| < 101, 0 \leq t \leq
1$ when $0 \leq m \leq l-1$;
\begin{equation}  \label{expest}
|\frac{\partial^m (y_k -v_k)}{\partial v_{i_1} \cdots
\partial v_{i_m}}| \leq  {\rm constant} \cdot  \varepsilon  \mbox{,}
\hspace{.5in} |\frac{\partial^m (x_k -p_k - v_k t)}{\partial v_{i_1}
\cdots \partial v_{i_m}}| \leq  {\rm constant} \cdot  \varepsilon.
\end{equation}
We will do this by induction on $m \geq 0$;
%And setting \( H_1(t)= \sup_{k} {\displaystyle \sum_{i=1}^4 \sup_{0 \leq s \leq t} | \frac{\partial y_i}{\partialv_k} (s,v)| } \),
 The case $m=0$ holds by (\ref{idon}).
 Suppose that (\ref{expest}) holds for $m=n-1 \leq l-2$, for some integer $n \geq
 1$. Then  $|\int_0^t L_n ds| \leq c_6  \varepsilon$ from induction hypothesis.
From the first equation in (\ref{dyv}), $|\frac{\partial^n
x_k}{\partial v_{i_1}  \cdots \partial v_{i_n}}(t, v)| \leq \sup_{s
\in [0,t]} |\frac{\partial^n y_k}{\partial v_{i_1}  \cdots
\partial v_{i_n}}(s,v)|$, as $0 \leq t \leq 1$. So, using (\ref{idon}),
from the second equation in (\ref{dyv}) we deduce;
\begin{equation}  \label{117a} \sum_{k=1}^{4}
\sup_{s \in [0,t]} |\frac{\partial^n (y_k -v_k)}{\partial v_{i_1}
\cdots \partial v_{i_n}}   (s, v)| \leq   c_7 \varepsilon
\sum_{k=1}^{4} \sup_{s \in [0,t]} |\frac{\partial^n y_{k} }{\partial
v_{i_1} \cdots
\partial v_{i_n}} (s, v)| + 4 c_6 \varepsilon.
\end{equation}

\noindent From this, for small  $\varepsilon$ we easily derive
(\ref{expest}) when $m=n$. Then,  for any $0 \leq m \leq l-1 $,
\begin{equation} \label{121}
|\frac{\partial^m (f_k(v)-(p_k +v_k) )}{\partial v_{i_1}  \cdots
\partial v_{i_m}} | \leq c_{8} \varepsilon  \ \ \mbox{ and } \ \
|\frac{\partial^m (h_k(x) -(x_k - p_k) )}{\partial x_{j_1} \cdots
\partial x_{j_m}}| \leq c_{9} \varepsilon.
\end{equation}
%\hspace{.1in}  & |\frac{\partial y_l}{\partial v_k} (t, v)- \delta_k^l| \leq c_7 \varepsilon \mbox{,} \hspace{.4in}
%|\frac{\partial x_l}{\partial v_k} (t, v)- \delta_k^l t| \leq c_7 \varepsilon t   \mbox{} \hspace{.3in}\nonumber \\
%|\frac{\partial f_l}{\partial v_k} (v)- \delta_k^l|  \leq  c_7 \varepsilon,  \hspace{.5in}
% |\frac{\partial h_l}{\partial x_k}- \delta_k^l| \leq c_{8} \varepsilon \mbox{,} \hspace{.4in} \end{align}
%\sum_{i=1}^4 \sup_{0 \leq s \leq t} | \frac{\partial y_i}{\partial v_k} (s)|,
where the second inequality is deduced from the first, using
(differentiation of)  $ \ h \circ f (v) = v$.  We let
$\lambda^{(0)}$ to be a small $\varepsilon$ value which makes all
above arguments hold. Choice of $C^{(0)}$ should be obvious. This
proves the lemma.
%\begin{equation} \sum_{k=1}^{4}  \frac{\partial h_i}{\partial x_k} \cdot \frac{\partial f_k}{\partial
% v_j}= \delta_{j}^{i}. \label{invfg1} \end{equation}
%\noindent For the $g_{0}$-unit vector $u=\frac{\partial}{\partialx_k}$ at $q \in B^g_{100}(p) - \{ p \}$.
%So we show $| u(r_g) - u(r_{g_0})| < \rm{constant} \cdot \lambda$ for $u= \frac{\partial}{\partial x_k }$.
%Here$C{(\varepsilon)}$'s may have different values.
%Applying  $\frac{\partial^k }{\partial v_{i_1} \cdots \partial v_{i_k}}$
\end{proof}
\begin{lemma} \label{lemdr}
There exist $\lambda^{(1)}$ with $\lambda^{(1)} < \lambda_1$ and a
positive constant $C^{(1)}$ such that  if $\|g- g_{0}
\|_{{C^l_{g_{0}}}} \leq \lambda \leq \lambda^{(1)}$ on
$B^{g_{0}}_{10^4}(p) \subset {\Bbb R}^n$, where $4 \leq l \leq 2
\kappa +4$, then for $j=0, 1, \cdots ,l-2$,
 \begin{equation} \label{lemdr1}
  \|(\nabla^{g_0})^j dr_g - (\nabla^{g_0})^j dr_{g_0} \|_{{C^0_{g_{0}}}}
< \frac{C^{(1)} \cdot {\lambda}}{r_{g_0}^j} \hspace{0.2in} \mbox{ on
} B^g_{100}(p) - \{ p \}.
\end{equation}
%, where  $C{(\varepsilon)}$ is a constant depending on $\varepsilon$ and
% ${\displaystyle \lim_{\varepsilon \rightarrow 0} C{(\varepsilon)} =0}$.
\end{lemma}

\begin{proof}
We choose $\lambda^{(1)}$ to be $\lambda^{(0)}$ of Lemma \ref{106}.
 We have $r_g (x)=\mid{\rm exp}_p^{-1} (x)|_{g|_p} =\sqrt{
\sum_{i, j=1}^4 g_{ij} (p) \mbox{ } h_i(x) h_j(x) }$ for $x$ with
$r_g (x)  < 200$.
 %Note that by Lemma \ref{lamb}, $ t |v|_{g|_p} $ equals $r_g(p, x(t,v))$ if the latter
 %is less than or equal to $200$.
We shall prove (\ref{lemdr1}) by induction on $j$. The case $j=0$
can be checked easily by applying (\ref{121}).
%\begin{align*}  |r_g(x) -r_{g_0}(x)|  &= \mid \sqrt{\sum_{i,j} g_{ij}(p) h_i h_j} - |x-p| \mid \mbox{ } \leq c_8 \varepsilon r_{g_0}, \\
% |\ \frac{\partial r_g}{\partial x_k } -  \frac{\partial r_{g_{0}}}{\partial x_k } |= &  |\frac{r_{g_{0}}  \sum_{i,j}g_{ij}(p) h_i
%(\frac{\partial}{\partial x_k} h_j) - r_g  (x_k-p_k)} {r_g r_{g_{0}} }  | \nonumber
%\leq \frac{ c_{9} \varepsilon r_{g_0}^2}{(1-c_{9}\varepsilon) r_{g_0}^2} < c_{10}\varepsilon . \end{align*}
Suppose that (\ref{lemdr1}) holds for $j=0, \cdots, i-1 \leq l-3$.
For the rest of this proof we write $\partial^k f$ to denote a
$k$-th order partial derivative of a function $f$ in $x$
coordinates, by an abuse of terminology. Then we may write for $\mu=
g$ or $g_0$
 \begin{equation} \label{33}
 \partial^{i+1} (r_{\mu}^2)  -2r_{\mu}(\partial^{i+1} r_{\mu})   =
\sum_{k=1}^{i} a_k (\partial^k r_{\mu}) (\partial^{i+1-k} r_{\mu})
\end{equation}
 where
  $a_k$'s are natural numbers. Here $\partial^{i+1}$ in the LHS is a specific partial
  derivative, e.g.
  $\frac{\partial^{i+1} r_{\mu}^2}{\partial x_{s_1} \cdots
\partial x_{s_{i+1}}}$ but  the index $k$ in the RHS should be understood as a
multi-index.
%When $j=0$ or $1$, this clearly holds. Assume that it holds for $j =0, 1, \cdots, i-1$.
%From (\ref{33}), $|r_{g_0}(\partial^{k+1} r_{g_0})|   =\sum_{s=1}^{k} \frac{a_s}{2} | (\partial^s r_{g_0})|
%\cdot |(\partial^{k+1-s} r_{g_0})| \leq \frac{c(i)}{r_{g_0}^{k}}$.
We have
\begin{equation} \label{zz}
\partial^{i+1} r_g - \partial^{i+1} r_{g_0} -\frac{\partial^{i+1}
({r_g}^{ 2})}{2r_g} + \frac{\partial^{i+1} (r_{g_0}^2)}{2r_{g_0}}  =
\sum_{k=1}^{i} \frac{a_k}{2} (\frac{\partial^k r_{g_0}
\partial^{i+1-k} r_{g_0}}{r_{g_0}} - \frac{\partial^k r_g \partial^{i+1-k}
r_g}{r_g}).
\end{equation}
 Using $(\ref{33})$ for $\mu= g_0$ and
$\partial^{i+1} r_{g_0}^2 =0$ for $i \geq 2$, one can derive, by
induction;
\begin{equation}  \label{35}
| \partial^{i+1} r_{g_0} | \leq \frac{c_{10}}{r_{g_0}^{i}} \mbox{
for } i =0, 1, \cdots, l-2.
\end{equation}
 Using this inequality and the induction
hypothesis, $\sum_{k=1}^{i} a_k |\frac{\partial^k r_{g_0}
\partial^{i+1-k} r_{g_0}}{r_{g_0}} - \frac{\partial^k r_g \partial^{i+1-k}
r_g}{r_g}| \leq \frac{ c_{11} \lambda}{r_{g_0}^{i}}$. We need to
show $|\frac{\partial^{i+1} (r_g^2)}{r_g} - \frac{\partial^{i+1}
(r_{g_0}^2)}{r_{g_0}}|\leq \frac{{\rm constant} \cdot
\lambda}{r_{g_0}^{i}}$.  One can check easily from (\ref{121}) that
  $r_g^2 =\sum_{i, j=1}^4 g_{ij}(p) \mbox{ } h_i(x)
h_j(x) $ satisfies $|\frac{\partial^{2} (r_g^2)}{r_g} -
\frac{\partial^{2} (r_{g_0}^2)}{r_{g_0}}|\leq \frac{ c_{12}
\lambda}{r_{g_0}}$.
 For $i \geq 2$,
 $|\partial^{i+1} (r_g^2)| \leq c_{13} \lambda$  from
(\ref{121}). As $\partial^{i+1} (r_{g_0}^2)=0$, we have
$|\frac{\partial^{i+1} (r_g^2)}{r_g} - \frac{\partial^{i+1}
(r_{g_0}^2)}{r_{g_0}}|= |\frac{\partial^{i+1} (r_g^2)}{r_g} | \leq
\frac{c_{13} \lambda}{r_g} \leq \frac{c_{14} \lambda}{r_{g_0}^{i}}$
 on $B^g_{100}(p) -\{p\}$. In sum, from (\ref{zz}) we get $|\partial^{i+1} r_g -
\partial^{i+1} r_{g_0}| \leq \frac{c_{15} \lambda}{r_{g_0}^{i}}$. This proves (\ref{lemdr1}) for some $C^{(1)}$ when $0 \leq  j
\leq l-2 \leq 2 \kappa +2$ and the Lemma.
\end{proof}

\noindent {\bf Proof of Proposition \ref{lemom};} In Lemma
\ref{lemdr}  we already proved the inequality (\ref{prop1}) for
$\omega_1$. From the hypothesis on $g$ we have $ \|J_g - J_0
\|_{{C^l_{g_{0}}}} \leq \lambda$, so the proposition holds for
$\omega_2= J_g \omega_1$. With these results for $\omega_1$ and
$\omega_2$, by choosing $\lambda_2$ small enough, the denominator of
$\omega_3$ is never zero on $B^g_{100} (p) - \{ p \}$. It is now a
straightforward computation to show (\ref{prop1}) for $\omega_3$ and
$\omega_4$.
%the following for $\omega_3=\frac{ \omega_3^0 -g(\omega_3^0,\omega_1)\omega_1 - g(\omega_3^0, \omega_2)\omega_2}{|    \omega_3^0
%-g(\omega_3^0, \omega_1)\omega_1 - g(\omega_3^0, \omega_2)\omega_2 |_g } $; $|\frac{\partial \omega_3}{\partial x_k} - \frac{\partial
%\omega^0_3}{\partial x_k}| < \frac{C{ \cdot \lambda}}{r_{g_0}}$, $\mbox{ }|\frac{\partial^2 \omega_3}{\partial x_l \partial x_k } -
%\frac{\partial^2 \omega^0_3}{\partial x_l \partial x_k }| < \frac{C{ \cdot \lambda}}{r_{g_0}^2}, $ from which the formula (\ref{prop1})follows.
This finishes the proof of Proposition \ref{lemom}. \hfill \quad

\bigskip
Now we shall derive the formula (\ref{sgt}) in section 4.
 For the metric $g= \sum_{i=1}^4
\omega_i \otimes \omega_i$ in section 4 we write $d\omega_i =
\sum_{j} \omega_{ij} \wedge \omega_j$ with $ \omega_{ij}
=-\omega_{ji}$  and similarly for $\tilde{g} :=g_{d,s}^{b, \epsilon,
c}$, $d \tilde{\omega}_i = \sum_{j} \tilde{\omega}_{ij} \wedge
\tilde{\omega}_j$ with $ \tilde{\omega}_{ij} =-\tilde{\omega}_{ji}$.
We set  $\omega_{ij}= \sum_{k} a_{ijk} \omega_{k}$ and
$\tilde{\omega}_{ij}= \sum_{k} \tilde{a}_{ijk} \tilde{\omega}_{k}$.
%Clearly, $a_{ijk} = -a_{jik}$ and $\tilde{a}_{ijk} =-\tilde{a}_{jik}$.
 Then,
\begin{equation} \label{comp}
 d\omega_{i}= \sum_{j<k} (a_{ikj} - a_{ijk}) \omega_j \wedge
\omega_k \ \  \mbox{  and  } \ \  d \tilde{\omega}_i = \sum_{j<k}
(\tilde{a}_{ikj} - \tilde{a}_{ijk}) \tilde{\omega}_j \wedge
\tilde{\omega}_k.
\end{equation}
 From (\ref{oms}),  $d\omega_1 = d(dr_g)= \sum_{j<k} (a_{1kj} - a_{1jk}) \omega_j \wedge \omega_k
 =0$,
so  $a_{1jk}= a_{1kj}$ for any $j, k$. We let $\alpha(r_g) =
\frac{1}{\sqrt{1+ Y(r_g)}}= {1 \over {\sqrt{ 1+ h_{\epsilon}^{b}
\cdot f_{d,s}(c-r_g) }}}$ so that $ \tilde{\omega}_1 = \alpha(r_g)
\omega_1$. Then $d \tilde{\omega}_1 = d (\alpha(r_g) dr_g)=0$, so
$\tilde{a}_{1jk}= \tilde{a}_{1kj}$. Comparing the two equations in
(\ref{comp}), one gets the following relations. Other ones can be
derived from these, using $a_{ijk}= - a_{jik}$ and $a_{1jk}=
a_{1kj}$.
\begin{align}
 \mbox{ For } i=2,3,4, \ \ \tilde{a}_{1i1} =&  0, \ \ \ \ a_{1i1}=0.  \  \nonumber \\
 \mbox{For } j=3,4,  \ \ \ \ \tilde{a}_{2j2} = & a_{2j2}, \ \ \tilde{a}_{j1j} = \frac{1}{\alpha}
 a_{j1j},  \ \ \
 \tilde{a}_{j2j} =  \alpha a_{j2j},  \  \nonumber\\
 \tilde{a}_{2j1} = \frac{1}{2}(1 - \frac{1}{\alpha^2}) a_{21j} \ & + \frac{1}{2}(1 + \frac{1}{\alpha^2})
 a_{2j1},
 \ \
 \tilde{a}_{21j} =\frac{1}{2}(1+ \frac{1}{\alpha^2}) a_{21j} + \frac{1}{2}(1 - \frac{1}{\alpha^2})
 a_{2j1}.
 \nonumber \\
 \tilde{a}_{212} =  \frac{\alpha^{\prime}}{\alpha^2} +  \frac{a_{212}}{\alpha},   \
 \ \ \ \ &  \tilde{a}_{343} =  a_{343},   \ \
\tilde{a}_{434} = a_{434}, \ \  \tilde{a}_{143} = \frac{1}{\alpha}
 a_{143},  \ \ \tilde{a}_{341} = \frac{1}{\alpha}
 a_{341} \nonumber \\ \  \tilde{a}_{342} =
\frac{1}{2\alpha} (a_{243}- a_{234}&) + \frac{\alpha}{2} (2a_{342}+
a_{234}-a_{243}),
\nonumber \\
 \tilde{a}_{234} = \frac{1}{2\alpha} (a_{234}- a_{243}&) + \frac{\alpha}{2} (a_{243}+
 a_{234}), \  \nonumber \\
 \tilde{a}_{243} = \frac{1}{2\alpha} (a_{243}- a_{234}&) + \frac{\alpha}{2} (a_{243}+ a_{234}) .
 \nonumber   \label{ata}
\end{align}

\noindent  We compute the curvature components $\tilde{R}_{ijij}$ of
$g_{d,s}^{b, \epsilon, c}$ from $ d \tilde{\omega}_{ij} - \sum_{k}
\tilde{\omega}_{ik} \wedge \tilde{\omega}_{kj} = \sum_{k<l}
\tilde{R}_{ijkl} \tilde{\omega}_k \wedge \tilde {\omega}_l$ and the
scalar curvatue. We set $d(\tilde{a}_{ijk}) = \sum_{s}
\tilde{a}_{ijk,s} \tilde{\omega}_s$.
 %Note that $R_{1j1j}$, $j=2,3,4$, is simpler than others because $d\tilde{\omega}_{1}=0$.
\begin{align*}
\tilde{R}_{1212} = &\tilde{a}_{122,1}   + \tilde{a}_{122}^2 +
\tilde{a}_{123}^2 + \tilde{a}_{124}^2 + 2
(\tilde{a}_{123}\tilde{a}_{321}+\tilde{a}_{124} \tilde{a}_{421}), \\
\tilde{R}_{1313} = & \tilde{a}_{133,1}   + \tilde{a}_{132}^2 +
\tilde{a}_{133}^2 + \tilde{a}_{134}^2  + 2
(\tilde{a}_{132}\tilde{a}_{231}+\tilde{a}_{134} \tilde{a}_{431}),\\
\tilde{R}_{1414} =  &\tilde{a}_{144,1}   + \tilde{a}_{142}^2 +
\tilde{a}_{143}^2 + \tilde{a}_{144}^2  + 2
(\tilde{a}_{142}\tilde{a}_{241}+\tilde{a}_{143} \tilde{a}_{341}), \\
\tilde{R}_{2323} = &\tilde{a}_{233,2}  - \tilde{a}_{232,3} +
\tilde{a}_{234}(\tilde{a}_{432}- \tilde{a}_{423}) +
\tilde{a}_{232}^2 + \tilde{a}_{233}^2 +
\tilde{a}_{122}\tilde{a}_{133} +\tilde{a}_{422}\tilde{a}_{433} \\
 &- \tilde{a}_{123}^2 + \tilde{a}_{243} \tilde{a}_{432},\\
\tilde{R}_{2424} =& \tilde{a}_{244,2}  - \tilde{a}_{242,4} +
\tilde{a}_{243}(\tilde{a}_{342}- \tilde{a}_{324}) +
\tilde{a}_{242}^2 + \tilde{a}_{244}^2 +
\tilde{a}_{122}\tilde{a}_{144} +\tilde{a}_{322}\tilde{a}_{344} \\
 &- \tilde{a}_{124}^2 + \tilde{a}_{234} \tilde{a}_{342},\\
\tilde{R}_{3434} =& \tilde{a}_{344,3}  - \tilde{a}_{343,4} +
\tilde{a}_{342}(\tilde{a}_{243}- \tilde{a}_{234}) +
\tilde{a}_{343}^2 +  \tilde{a}_{344}^2 +
\tilde{a}_{133}\tilde{a}_{144} +\tilde{a}_{233}\tilde{a}_{244} \\
 &- \tilde{a}_{134}^2 + \tilde{a}_{324} \tilde{a}_{243},
\end{align*}
\begin{align*}
 \frac{s(g_{d,s}^{b, \epsilon, c})}{2}
         = & \tilde{R}_{2112} + \tilde{R}_{3113} + \tilde{R}_{4114} + \tilde{R}_{3223} + \tilde{R}_{4224} + \tilde{R}_{4334} \nonumber \\
         = & -\sum_{i=2}^4 \tilde{a}_{1ii,1} - \sum_{2 \leq i \leq j \leq 4} \tilde{a}_{1ij}^2
          + \sum_{i=3}^4 (\tilde{a}_{2i2,i} - \tilde{a}_{2ii,2}) + \tilde{a}_{343,4}  -\tilde{a}_{344,3}\nonumber \\
           & +\tilde{a}_{423}\tilde{a}_{342}
           + \tilde{a}_{234}\tilde{a}_{342} + \tilde{a}_{234}\tilde{a}_{423} - \sum_{i=3}^4(\tilde{a}_{2i2}^2 + \tilde{a}_{2ii}^2  )-\tilde{a}_{343}^2 - \tilde{a}_{344}^2
            \nonumber \\
            &  - \sum_{i=3}^4 \tilde{a}_{122}\tilde{a}_{1ii}  - \tilde{a}_{422}\tilde{a}_{433}  -\tilde{a}_{322}\tilde{a}_{344}
             -\tilde{a}_{133}\tilde{a}_{144} -
             \tilde{a}_{233}\tilde{a}_{244}\\
                       = & \ \frac{\alpha^{''}}{\alpha^3}- 3\frac{(\alpha^{'})^{2}}{\alpha^4}+
           \frac{\alpha^{'}}{\alpha^3} (3a_{122} + 2 a_{133} + 2 a_{144})- \frac{1}{\alpha^2} \sum_{i=2}^{4} (a_{1ii,1}
           +a_{1ii}^{2}) \nonumber  \\
          & - \sum_{i=3}^{4} \frac{1}{4} \{ (1 + \frac{1}{\alpha^2}) a_{21i}  +  (1 - \frac{1}{\alpha^2}) a_{2i1}\}^2
          - \frac{1}{\alpha^2} a_{134}^2 \\
          & + \sum_{i=3}^{4} ( a_{2i2,i}-\alpha^2 a_{2ii,2} )   + a_{343,4} +  a_{434,3}
          +(\frac{1}{2}-\frac{1}{4 \alpha^2} )(a_{243} - a_{234} )^2
          \\
           &+a_{342} (a_{234} - a_{243} )
          - \frac{\alpha^2}{4}(a_{243} + a_{234} )^2 - \sum_{i=3}^{4} ( \alpha^2 a_{2ii}^2 +a_{2i2}^2)-a_{343}^2
\\
          & - a_{344}^2 -\sum_{i=3}^{4} \frac{a_{122}}{\alpha^2} a_{1ii}
           +  a_{242} a_{433} +  a_{232} a_{344} - \frac{1}{\alpha^2} a_{133} a_{144}
          -  \alpha^2 a_{233} a_{244}.
\end{align*}

\noindent Putting $ \frac{1}{\alpha^2} = 1 + Y$,
$\frac{\alpha^{'}}{\alpha^3}=-\frac{1}{2}Y^{'}$ and $
\frac{\alpha^{''}}{\alpha^3}- 3\frac{(\alpha^{'})^{2}}{\alpha^4}=
-\frac{1}{2} Y^{''}$ into the above, we get the formula (\ref{sgt}).

\end{document}